\documentclass[11pt]{article}
\usepackage{amssymb,amsmath,amsthm,cite}
\newcommand{\R}{\mathbb{R}}
\newcommand{\Z}{\mathbb{Z}}
\newcommand{\N}{\mathbb{N}}
\def\qed{\hfill $\Box$ \smallskip}
\def\x#1{{\rm (\ref{#1})}}

\allowdisplaybreaks[4]

\begin{document}
\title{\LARGE\bf{ Existence of ground state solutions of Nehari-Pankov type to Schr\"odinger systems}
 \footnote{This work is partially supported by the NNSF (No: 11571370, 11471137) of China. This paper has been accepted for publication in SCIENCE CHINA Mathematics.}}
\date{}
 \author{TANG XianHua$^{1}$ and LIN XiaoYan$^{2}$  \\
        {\small $^{1}$ School of Mathematics and Statistics, Central South University,}\\
        {\small Changsha, Hunan 410083, P.R.China }\\
        {\small E-mail: tangxh@mail.csu.edu.cn}\\
        {\small $^{2}$ Department of Mathematics, Huaihua University,}\\
        {\small Huaihua, Hunan 418008, P.R.China}\\
        {\small E-mail: xiaoyanlin98@hotmail.com}\\
}
\maketitle
\begin{center}
\begin{minipage}{13cm}
\par
\small  {\bf Abstract:} This paper is dedicated to studying the following elliptic system of Hamiltonian type:
$$
   \left\{
   \begin{array}{ll}
    -\varepsilon^2\triangle u+u+V(x)v=Q(x)F_{v}(u, v), \ \ \ \  x\in {\R}^{N},\\
    -\varepsilon^2\triangle v+v+V(x)u=Q(x)F_{u}(u, v), \ \ \ \  x\in {\R}^{N},\\
    |u(x)|+|v(x)| \rightarrow 0, \ \ \mbox{as} \ |x|\rightarrow \infty,
  \end{array}
    \right.
 $$
 where $N\ge 3$, $V, Q\in \mathcal{C}(\R^N, \R)$, $V(x)$ is allowed to be sign-changing and $\inf Q > 0$, and
 $F\in \mathcal{C}^1(\R^2, \R)$ is superquadratic at both $0$ and infinity but subcritical. Instead of the reduction
 approach used in [Calc Var PDE, 2014, 51: 725-760], we develop a more direct approach -- non-Nehari manifold approach to obtain stronger conclusions but under weaker assumptions than these in [Calc Var PDE, 2014, 51: 725-760]. We can find an $\varepsilon_0>0$ which is determined by terms of $N, V, Q$ and $F$, then we prove the
 existence of a ground state solution of Nehari-Pankov type to the coupled system for all $\varepsilon\in (0, \varepsilon_0]$.


\end{minipage}
\end{center}

 {\section{Introduction}}
 \setcounter{equation}{0}

 \par
   In this paper we study standing waves for the following system of time-dependent nonlinear Schr\"odinger equations
 \begin{equation} \label{HS1}
   \left\{
   \begin{array}{ll}
    i\hbar \frac{\partial \phi_1}{\partial t}+\frac{\hbar^2}{2m}\triangle \phi_1+\phi_1+f(x, \phi)\phi_2=0,\\
    i\hbar \frac{\partial \phi_2}{\partial t}+\frac{\hbar^2}{2m}\triangle \phi_2+\phi_2+f(x, \phi)\phi_1=0,
   \end{array}
    \right.
 \end{equation}
 where $m$ is the mass of a particle, $\hbar$ is the Planck constant, $\phi = (\phi_1, \phi_2)$,
 $\phi_1(t, x)$ and $\phi_1(t, x)$ are the complex valued envelope functions. Suppose that $f (x, e^{i\theta}\phi) =
 f(x, \phi)$ for $\theta \in [0, 2\pi]$. We will look for standing waves of the form
 $$
   \phi_1 (t, x) = e^{i\omega t}u(x), \ \ \ \ \phi_2 (t, x) = e^{i\omega t}v(x),
 $$
 which propagate without changing their shape and thus have a soliton-like behavior.  System \x{HS1} arises quite
 naturally in nonlinear optics and Bose-Einstein condensates (see \cite{AC, LW, SS} and the references therein). In
 general, the above coupled nonlinear Schr\"odinger system leads to the elliptic system of
 Hamiltonian form
 \begin{equation} \label{HS2}
   \left\{
   \begin{array}{ll}
    -\varepsilon^2\triangle u+u=H_{v}(x, u, v), \ \ \ \  x\in {\R}^{N},\\
    -\varepsilon^2\triangle v+v=H_{u}(x, u, v), \ \ \ \  x\in {\R}^{N},\\
    |u(x)|+|v(x)| \rightarrow 0, \ \ \mbox{as} \ |x|\rightarrow \infty,
   \end{array}
    \right.
 \end{equation}
 where $N \ge 3$, $H \in \mathcal{C}^1(\R^N \times \R^2, \R)$ and $\varepsilon > 0$ is a small parameter. The study of the systems similar to \x{HS2} has only begun quite recently. When $\varepsilon = 1$, it was considered recently
 in some works \cite{Bd, dF, LY, LTZ, ZCZ, ZQZ, ZTZ, ZZD, ZD, ZZ}. For a similar problem on a bounded domain we
 refer the reader to \cite{BR, Cd, dF, HV} and the references therein. For a survey on this direction see \cite{de, Ru}.

 \par
   For the case $\varepsilon > 0$ is a small parameter, there are some recent works considering the
 existence of solutions;
 see for instance \cite{AS, DL, SS} and the references therein. In contrast with the case $\varepsilon = 1$,
 except for the difficulties that the lack of the compactness of the Sobolev embedding and the energy functional
 is strongly indefinite, no uniqueness results seem to be known for the ``limit problem" and this is in some cases
 a crucial assumption in the single equation case. So asymptotic analysis of solutions with respect to small
 $\varepsilon > 0$ has been very recently performed; see for example \cite{AS, AY, DLZ, DL, PR, RS, RT, RY, SS}
 and their references. Except for \cite{DLZ, DL, SS}, most of the above works considered the
 case that $H(x, u, v) = F(u) + G(v)$. In particular, in \cite{AY}, the authors obtained the existence
 of positive solutions which concentrate on the boundary of $\Omega$ for an elliptic system with zero
 Neumann boundary condition on a bounded domain $\Omega$ (see also \cite{PR}). In \cite{RS}, Ramos and
 Soares considered the following problem
 \begin{equation} \label{HS3}
   \left\{
   \begin{array}{ll}
    -\varepsilon^2\triangle u+V(x)u=g(v), \ \ \ \  x\in \Omega,\\
    -\varepsilon^2\triangle v+V(x)v=f(u), \ \ \ \  x\in \Omega,\\
    u, v\in H^1(\R^N),
   \end{array}
    \right.
 \end{equation}
 where $\Omega$ is a domain of $\R^N$, $V\in \mathcal{C}(\R^N, \R)$ satisfies $0<V(0)=\min V(x)<
 \liminf_{|x|\to \infty}V(x)\in (0, \infty]$, $f(u)$ and $g(v)$ are power type functions, superlinear but subcritical
 at infinity. The authors established the existence of positive solutions
 which concentrate, as $\varepsilon \rightarrow 0$, at a prescribed finite number of local minimum points (possibly degenerate) of the
 potential $V$. Different from those discussed in \cite{AS, RS, RT}, Ding, Lee and Zhao \cite{DLZ} dealt with existence
 and concentration phenomena of the ground state solutions to the following subcritical problem
 \begin{equation} \label{HS4}
   \left\{
   \begin{array}{ll}
    -\varepsilon^2\triangle u+u+V(x)v=Q(x)g(|z|)v, \ \ \ \  x\in {\R}^{N},\\
    -\varepsilon^2\triangle v+v+V(x)u=Q(x)g(|z|)u, \ \ \ \  x\in {\R}^{N},\\
    |u(x)|+|v(x)| \rightarrow 0, \ \ \mbox{as} \ |x|\rightarrow \infty,
   \end{array}
    \right.
 \end{equation}
 where $z:=(u, v)$, $V, Q\in \mathcal{C}^1(\R^N, \R)$ and $g\in \mathcal{C}^1(\R^{+}, \R^{+})$. Since the energy functional $\Phi_{\varepsilon}$
 associated with system \x{HS4} is strongly indefinite,
 to overcome this difficulty, as in \cite{AN} (see also \cite{RS} and \cite{RT}), the authors constructed a
 reduced functional $\mathcal{R}_{\varepsilon}$ whose critical points are in one to one to critical points
 of $\Phi_{\varepsilon}$, which was first proposed in \cite{AN}. With the help of the Nehari manifold
 of $\Phi_{\varepsilon}$, an important information of the least energy $c_{\varepsilon}$ was obtained.
 By estimating the asymptotic behavior of $c_{\varepsilon}$ as $\varepsilon \to 0$, they proved $c_{\varepsilon}$
 is attained for sufficiently small $\varepsilon>0$. In order to state their results, some notations and assumptions
 are required. Set
 $$
   V_{\min}:=\min_{x\in \R^N} V(x), \ \ \ \ V_{\max}:=\max_{x\in \R^N} V(x), \ \ \ \ \mathcal{V}:=\{x\in \R^N : V(x)= V_{\min}\};
 $$
 $$
   Q_{\min}:=\min_{x\in \R^N} Q(x), \ \ \ \ Q_{\max}:=\max_{x\in \R^N} Q(x), \ \ \ \ \mathcal{Q}:=\{x\in \R^N : Q(x)= Q_{\max}\};
 $$
 $$
   \ \ \ \ V_{\infty}:=\liminf_{|x|\to\infty}V(x),  \ \ \ \ Q_{\infty}:=\limsup_{|x|\to\infty}Q(x);
 $$
 $$
   \mathcal{A}_{v}:=\{x\in \mathcal{V} : Q(x)=Q(x_{v})\}\cup \{x\not\in \mathcal{V} : Q(x)> Q(x_{v})\}
 $$
 and
 $$
   \mathcal{A}_{q}:=\{x\in \mathcal{Q} : V(x)=V(x_{q})\}\cup \{x\not\in \mathcal{W} : V(x) < V(x_{q})\}.
 $$

 \par
   Furthermore, the following assumptions are required:

 \vskip2mm
 \noindent
 \begin{itemize}
 \item[(A0)]  $V, Q\in \mathcal{C}^1(\R^N, \R)$, $\|V\|_{\infty}< 1$, $0<Q_{\min}\le Q_{\max}<\infty$;

 \item[(A1)]  $V_{\min}< V_{\infty}$, and there exist $x_{v}\in \mathcal{V}$ and $R>0$ such that
 $$
   Q(x_{v})=\max_{y\in \mathcal{V}}Q(y)\ge Q(x), \ \ \ \ \forall \  |x|\ge R;
 $$

 \item[(A2)]  $Q_{\max} > Q_{\infty}$, and there exist $x_{q}\in \mathcal{Q}$ and $R>0$ such that
 $$
   V(x_{q})=\min_{y\in \mathcal{Q}}V(y)\le V(x), \ \ \ \ \forall \  |x|\ge R;
 $$

 \item[(G1)] $g\in \mathcal{C}^1(\R^{+}, \R^{+}), g(0)=0$ and $g'(s)\ge 0$ for $s>0$, where $\R^{+}:=[0, \infty)$;

 \item[(G2)] there exist $C>0$ and $p\in (2, 2^*)$ such that $|g(s)|\le C(1+s^{p-2})$ for all $s\in \R^{+}$;

 \item[(G3)] there exists $\mu>2$ such that $ g(s)s^2\ge \mu \int_{0}^{s}g(t)t\mathrm{d}t>0$ if $s> 0$.
 \end{itemize}

 \vskip4mm
 \par
    In \cite{DLZ}, they proved the following theorem.

 \vskip4mm
 \par\noindent
 {\bf Theorem 1.1.} (\cite[Theorem 1]{DLZ})\ \  {\it Let } (A0), (G1), (G2) {\it and}
 (G3) {\it be satisfied. Suppose that} (A1) {\it or} (A2) {\it is satisfied. Then for sufficiently small $\varepsilon > 0$,
 \x{HS4} has a least energy solution $\hat{z}_{\varepsilon}=(\hat{u}_{\varepsilon}, \hat{v}_{\varepsilon})$. }

 \vskip4mm
 \par
    Theorem 1.1 is very interesting. In its proof, many new tricks were used to overcome the difficulties caused by
 the strong indefinity of the energy functional $\Phi_{\varepsilon}$ associated with system \x{HS4}. We point that
 the regularity assumptions $V, Q\in \mathcal{C}^1$ and $g\in \mathcal{C}^1$ are very crucial in \cite{DLZ}, which
 seem to be necessary when the reduction method is used.
 Motivated by the works \cite{DLZ}, in this paper, we further study the existence of the
 ground state solutions of Nehari-Pankov type to the following more general problem
 \begin{equation} \label{PHS}
   \left\{
   \begin{array}{ll}
    -\varepsilon^2\triangle u+u+V(x)v=Q(x)F_{v}(u, v), \ \ \ \  x\in {\R}^{N},\\
    -\varepsilon^2\triangle v+v+V(x)u=Q(x)F_{u}(u, v), \ \ \ \  x\in {\R}^{N},\\
    |u(x)|+|v(x)| \rightarrow 0, \ \ \mbox{as} \ |x|\rightarrow \infty,
   \end{array}
    \right.
 \end{equation}
 where $V, Q\in \mathcal{C}(\R^N, \R)$ and $F\in \mathcal{C}^1(\R^2, \R)$. Instead of the reduction method used in \cite{DLZ}, we will use a
 more direct approach -- non-Nehari manifold approach which was first proposed
 in \cite{Ta1} for a single Schr\"odinger equation (see also \cite{Ta2, Ta3}). The main idea of this approach
 is to construct a minimizing Cerami sequence for the energy functional outside Nehari-Pankov manifold by using the
 diagonal method, which is completely different from that of Szulkin and Weth \cite{SW}.  This approach is
 valid when finding a ground state solution of Nehari-Pankov type.

 \par
    We will obtain stronger conclusions on existence of the ground state solutions of Nehari-Pankov type to \x{PHS}
 for small $\varepsilon>0$ but under weaker assumptions than these in \cite{DLZ}. Roughly speaking, we can find an
 $\varepsilon_0>0$ which is determined by terms of $N, V, Q$ and $F$, then we prove the existence of a ground state
 solutions of Nehari-Pankov type to \x{PHS} for all $\varepsilon\in (0, \varepsilon_0]$. In particular, we only need
 $V, Q\in \mathcal{C}(\R^N, \R)$ and $F\in \mathcal{C}^1(\R^2, \R)$. To the best of our knowledge, there seems to be no similar
 results in literature.

 \par
    To state our theorems accurately, we set

 \begin{equation}\label{ND0}
  \mathcal{ND}_{0}=\left\{h\in \mathcal{C}(\R^{+}, \R^{+}): \left\{\begin{array} {ll} h(0)=0 \ \mbox {and} \ h(s)\ \mbox{is nondecreasing on}\ \R^{+},\\
    \mbox{there exist constants}\ p \in (2, 2^*) \ \mbox{and}\ c_0>0 \ \mbox{such that}\\
     |h(s)|\le c_0(1+|s|^{p-2}),\ \ \forall \ s\ge 0.
  \end{array}\right.\right\}
 \end{equation}
 \vskip2mm
 \noindent
  Furthermore, we make the following assumptions.

 \vskip2mm
 \noindent
 \begin{itemize}
 \item[(V0)]  $V, Q\in \mathcal{C}(\R^N, \R)$, $\|V\|_{\infty}\le \frac{2\eta ab}{a^2+b^2}$ and $0<Q_{\min}\le Q_{\max}<\infty$,
 where $a, b>0$ and $\eta\in (0, 1)$;

 \item[(V1)]  $V_{\min}< V_{\infty}$, and there exist $x_{v}\in \mathcal{V}$ and $R>0$ such that
 $$
   Q(x_{v})\ge Q(x), \ \ \ \ \forall \  |x|\ge R;
 $$

 \item[(V2)]  $Q_{\max} > Q_{\infty}$, and there exist $x_{q}\in \mathcal{Q}$ and $R>0$ such that
 $$
   V(x_{q})\le V(x), \ \ \ \ \forall \  |x|\ge R;
 $$

 \item[(F1)] there exist $g_i, h_j\in \mathcal{ND}_0$, $\alpha_i, \beta_i, \alpha_j', \beta_j'\in \R$ with
 $\alpha_i^2+\beta_i^2\ne 0$ and $\alpha_j'>{\beta_j'}^2$, $i=1, 2, \ldots, k; j=1, 2, \ldots, l$, such that
 $$
   F(u, v)=\sum_{i=1}^{k}\int_{0}^{|\alpha_iu+\beta_iv|}g_i(s)s\mathrm{d}s
     +\sum_{j=1}^{l}\int_{0}^{\sqrt{u^2+2\beta_j'uv+\alpha_j'v^2}}h_j(s)s\mathrm{d}s;
 $$

 \item[(F2)] $\lim_{|au+bv|\to \infty}\frac{|F(u, v)|}{|au+bv|^2}=\infty$;

 \item[(F3)] there exist $\mathcal{C}_0>0, T_0>0$, $\mu>2$  and $F_0\in \mathcal{C}(\R^2, \R)$ with $F_0(u, v)>0 $ if $au+bv\ne 0$,
 such that
 $$
   F(tz)\ge \mathcal{C}_0t^{\mu}F_0(z), \ \ \ \ \forall \ z\in \R^2, \ t\ge T_0.
 $$
 \end{itemize}

 \vskip4mm
 \par\noindent
 {\bf Remark 1.2.}\ \  It is clear that (i) \ (F3) is weaker than (AR)-condition: there exists $\mu>2$
 such that $F_{z}(z)\cdot z\ge \mu F(z)>0$ for $z\ne 0$; \  (ii) \ (F2) is also weaker than the common
 super-quadratic condition (SQ): $\lim_{|z|\to \infty}\frac{|F(z)|}{|z|^2}=\infty$;
 \  (iii) \ Let $F(z)=\int_{0}^{|z|}g(s)s\mathrm{d}s$. Then (G1)-(G3) imply (F1)- (F3);
 \ (iv) \ Let $Q(x)\equiv 1$ and $F(z)=\int_{0}^{|u|}g(s)s\mathrm{d}s+\int_{0}^{|v|}h(s)s\mathrm{d}s$, then problem \x{PHS}
 reduces to \x{HS3}, which was studied in \cite{RS}. Moreover, the assumptions in \cite[(H)]{RS} also imply
 (F1)-(F3).

 \vskip2mm
 \par
    Before presenting our results, we give three nonlinear examples to illustrate the above assumptions.

 \vskip4mm
 \par\noindent
 {\bf Example 1.3.}\ \ Let $F(z)=|au+bv|^{\mu}$, where $\mu \in (2, 2^*)$ and $a, b>0$ with $\|V\|_{\infty}<\frac{2ab}{a^2+b^2}$.
  It is easy to see that $F(z)$ satisfies (F1)-(F3) with $F_0=F$, but not (AR).

 \vskip4mm
 \par\noindent
 {\bf Example 1.4.}\ \ Let $F(z)=|au+bv|^{\mu}+(u^2+uv+v^2)\ln (1+u^2+uv+v^2)$, where $\mu \in (2, 2^*)$ and $a, b>0$ with $\|V\|_{\infty}<\frac{2ab}{a^2+b^2}$. It is easy to see that $F(z)$ satisfies (F1)-(F3) with
 $F_0(u, v)=|au+bv|^{\mu}$, but not (AR).

 \vskip4mm
 \par\noindent
 {\bf Example 1.5.}\ \ Let $F(z)=|2u+v|^{\mu}+|u+2v|^{\mu}$, where $\mu \in (2, 2^*)$.
  It is easy to see that $F(z)$ satisfies (F1)-(F3) with $a=b=1$ and $F_0=F$.

 \vskip4mm
 \par
   Let $x_{m}=x_{v}$ if (V1) holds, or $x_{m}=x_{q}$ if (V2) holds. Replacing $u(\varepsilon x+x_{m})$
 and $v(\varepsilon x+x_{m})$ by $u(x)$ and $v(x)$, respetively, then system \x{PHS} is equivalent to
 \begin{equation} \label{HS5}
   \left\{
   \begin{array}{ll}
    -\triangle u+ u+V(\varepsilon x+x_{m})v=Q(\varepsilon x+x_{m})F_{v}(u, v), \ \ \ \  x\in {\R}^{N},\\
    -\triangle v+ v+V(\varepsilon x+x_{m})u=Q(\varepsilon x+x_{m})F_{u}(u, v), \ \ \ \  x\in {\R}^{N},\\
    |u(x)|+|v(x)| \rightarrow 0, \ \ \mbox{as} \ |x|\rightarrow \infty.
   \end{array}
    \right.
 \end{equation}

 \vskip4mm
 \par
   Let $E=H^1(\R^N)\times H^1(\R^N)$. Then $E$ is a Hilbert space with the standard inner product
 $$
   (z_1, z_2)_{H^1(\R^N)}=(u_1, u_2)_{H^1(\R^N)}+(v_1, v_2)_{H^1(\R^N)}, \ \ \ \ \forall \ z_i=(u_i, v_i)\in E, \ \ i=1, 2,
 $$
 and the corresponding norm
 $$
   \|z\|_{H^1(\R^N)}=\left(\|u\|^2_{H^1(\R^N)}+\|v\|^2_{H^1(\R^N)}\right)^{1/2}, \ \ \ \ \forall \ z=(u, v)\in E.
 $$
 Let $E=E^{-}\oplus E^{+}$ be an orthogonal decomposition, see Section 2. Define a functional
 \begin{equation}\label{Phe}
   \Phi_{\varepsilon}(z)=\int_{{\R}^N}\left[\nabla u\cdot \nabla v+uv+\frac{1}{2}V(\varepsilon x+x_{m})|z|^2\right]\mathrm{d}x
     -\int_{{\R}^N}Q(\varepsilon x+x_{m})F(z)\mathrm{d}x
 \end{equation}
 for all $z=(u, v)\in E$. Under assumptions (V0), (F1) and (F2), $\Phi_{\varepsilon}\in \mathcal{C}^1(E, \R)$ and
 \begin{eqnarray}\label{Phde}
   \langle \Phi_{\varepsilon}'(z), \varphi \rangle
    &  =  & \int_{{\R}^N}\left[\nabla u\cdot\nabla \psi+\nabla v\cdot\nabla \phi+(u\psi+v\phi)
               +V(\varepsilon x+x_{m})z\cdot \varphi\right]\mathrm{d}x\nonumber\\
    &     & \ \ \ \ -\int_{{\R}^N}Q(\varepsilon x+x_{m})F_z(z)\cdot \varphi\mathrm{d}x, \ \ \ \ \forall \ z=(u, v),\ \varphi=(\phi, \psi)\in E.
 \end{eqnarray}
 Let
 \begin{equation}\label{Nee-}
   \mathcal{N}_{\varepsilon}^{-}  = \left\{z\in E\setminus E^{-} : \langle \Phi_{\varepsilon}'(z), z \rangle
     =\langle \Phi_{\varepsilon}'(z), \zeta \rangle=0, \ \ \ \ \forall \ \zeta\in E^{-} \right\}.
 \end{equation}
 $\mathcal{N}_{\varepsilon}^{-}$ was first introduced by Pankov \cite{Pa}, which is a subset of the Nehari manifold
 \begin{equation}\label{Nee}
   \mathcal{N}_{\varepsilon}  = \left\{z\in E\setminus \{0\} : \langle \Phi_{\varepsilon}'(z), z \rangle=0 \right\}.
 \end{equation}

 \vskip4mm
 \par
   We are now in a position to state the first main result of this paper.

 \vskip4mm
 \par\noindent
 {\bf Theorem 1.6.}\ \ {\it Assume that $V$, $Q$ and $F$ satisfy} (V0), (V1) {\it and} (F1)-(F3). {\it Then there exists
 an $\varepsilon_0>0$ such problem \x{PHS} has a nontrivial solution
 $\hat{z}_{\varepsilon}=(\hat{u}_{\varepsilon}, \hat{v}_{\varepsilon})\in \mathcal{N}_{\varepsilon}^{-}$ with $\Phi_{\varepsilon}(z_{\varepsilon})=\inf_{\mathcal{N}_{\varepsilon}^{-}}\Phi_{\varepsilon}>0$ for
 $\varepsilon\in (0, \varepsilon_0]$, where $z_{\varepsilon}(x)=\hat{z}_{\varepsilon}(\varepsilon x+x_{v})$. }
 {\it If} (V1) {\it is replaced by} (V2), {\it then the above
 conclusion remains true by replacing $x_{v}$ with $x_{q}$.}

 \vskip4mm
 \par
   The ``limit problem" associated to \x{HS5} is an autonomous system
 \begin{equation} \label{HS07}
   \left\{
   \begin{array}{ll}
    -\triangle u+ u+V(x_{m})v=Q(x_{m})F_{v}(u, v), \ \ \ \  x\in {\R}^{N},\\
    -\triangle v+ v+V(x_{m})u=Q(x_{m})F_{u}(u, v), \ \ \ \  x\in {\R}^{N}.
   \end{array}
    \right.
 \end{equation}
 We will prove that the least energy $c_{\varepsilon}:=\inf_{\mathcal{N}_{\varepsilon}^{-}}\Phi_{\varepsilon}$ is attained for $\varepsilon\in (0, \varepsilon_0]$ by comparing with $c_{\varepsilon}$ and the least energy $c_0$
 associated with ``limit problem" \x{HS07}. Therefore, it is very crucial if $c_0$ can be attained, i.e. if \x{HS07} has a solution at which $\Phi_0$ has the least energy $c_0$ on $\mathcal{N}_{0}^{-}$. Prior to this, we consider the following more general periodic system
 \begin{equation} \label{HS}
   \left\{
   \begin{array}{ll}
    -\triangle u+V_1(x)u+V_2(x)v=W_{v}(x, u, v), \ \ \ \  x\in {\R}^{N},\\
    -\triangle v+V_1(x)v+V_2(x)u=W_{u}(x, u, v), \ \ \ \  x\in {\R}^{N},\\
    u, v\in H^{1}({\R}^{N}),
   \end{array}
    \right.
 \end{equation}
 where $N\ge 3$,  $V_1, V_2:{\R}^{N}\rightarrow \R$ and $W:{\R}^{N}\times\R^2\rightarrow \R$.
   More precisely, we make the following assumptions.

 \vskip2mm
 \noindent
 \begin{itemize}
  \item[(V0$'$)]  $V_1, V_2\in \mathcal{C}(\R^N)$ and satisfy
 \begin{equation}\label{V0}
  |V_2(x)|\le \frac{2\eta ab}{a^2+b^2}V_1(x), \ \ \ \ 0<\inf_{x\in \R^N}V_1(x)\le \sup_{x\in \R^N}V_1(x)<\infty,
 \end{equation}
 where $a, b>0$ and $\eta\in (0, 1)$;

 \item[(V1$'$)]   $V_1(x)$ and $V_2(x)$ are 1-periodic in each of $x_1, x_2, \ldots, x_N$;

 \item[(W0)] $W\in \mathcal{C}({\R}^{N}\times\R^2, \R^{+})$, $W(x, z)$ is continuously differentiable on $z\in \R^2$ for every
 $x\in \R^N$, and there exist constants $p \in (2, 2^*)$ and $C_0>0$ such that
 \begin{equation}\label{1.4}
   |W_{z}(x, z)|\le C_0\left(1+|z|^{p-1}\right),  \ \ \ \ \forall \ (x, z)\in \R^N\times \R^2;
 \end{equation}

 \item[(W1)] $W_{z}(x, z)=o(|z|)$, as $|z|\to 0$, uniformly in $x\in \R^N$;

 \item[(W2)] $\lim_{|au+bv|\to \infty}\frac{|W(x, u, v)|}{|au+bv|^2}=\infty$,\ a.e. $ x\in \R^N$;

 \item[(W2$'$)] $\lim_{|au+bv|\to \infty}\frac{|W(x, u, v)|}{|au+bv|^2}=\infty$, uniformly in  $x\in \R^N$;

 \item[(W3)] $W(x, z)$ is 1-periodic in each of $x_1, x_2, \ldots, x_N$;

 \item[(W4)] for all $\theta\ge 0, \ z, \zeta\in \R^2$, there holds
 $$
   \frac{1-\theta^2}{2}\nabla W(x, z)\cdot z-\theta\nabla W(x, z)\cdot \zeta+W(x, \theta z+\zeta)-W(x, z)\ge 0.
 $$
 \end{itemize}

 \vskip4mm
 \par
     Observe that, the natural functional associated with \x{HS} is given by
 \begin{equation}\label{Ph}
   \Phi(z)=\int_{{\R}^N}\left[\nabla u\cdot \nabla v+V_1(x)uv+\frac{1}{2}V_2(x)|z|^2\right]\mathrm{d}x-\int_{{\R}^N}W(x, z)\mathrm{d}x,
 \end{equation}
 for all $z=(u, v)\in E$. Furthermore, under assumptions (V0$'$), (W0) and (W1), $\Phi\in \mathcal{C}^1(E, \R)$ and
 \begin{eqnarray}\label{Phd}
   \langle \Phi'(z), \varphi \rangle
    &  =  & \int_{{\R}^N}\left[\nabla u\cdot\nabla \psi+\nabla v\cdot\nabla \phi+V_1(x)(u\psi+v\phi)
               +V_2(x)z\cdot \varphi\right]\mathrm{d}x\nonumber\\
    &     & \ \ \ \ -\int_{{\R}^N}W_z(x, z)\cdot \varphi\mathrm{d}x, \ \ \ \ \forall \ z=(u, v),\ \varphi=(\phi, \psi)\in E.
 \end{eqnarray}
 Let
 \begin{equation}\label{Ne-}
   \mathcal{N}^{-}  = \left\{z\in E\setminus E^{-} : \langle \Phi'(z), z \rangle=\langle \Phi'(z), \zeta \rangle=0,
       \ \forall \ \zeta\in E^{-} \right\}.
 \end{equation}

 \vskip2mm
 \noindent

 \vskip4mm
 \par
    For system \x{HS}, we obtain the following existence theorem on the ground state solutions of Nehari-Pankov
 type.

 \vskip4mm
 \par\noindent
 {\bf Theorem 1.8.}\ \ {\it Assume that $V$ and $W$ satisfy} (V0$'$), (V1$'$) {\it and} (W0)-(W4). {\it Then problem \x{HS} has a
 solution $z^*\in \mathcal{N}^{-}$ such that $\Phi(z^*)=\inf_{\mathcal{N}^{-}}\Phi>0$.}

 \vskip4mm
 \par
    However, it is not easy to check assumption (W4). Next, we give several classes functions satisfying (W4).
 Prior to this, we define one set as follows:

 \begin{equation}\label{ND}
  \mathcal{ND}=\left\{h\in \mathcal{C}(\R^{N}\times \R^{+}, \R^{+}): \left\{\begin{array} {ll} h(x, t)\ \mbox{is 1-periodic in each of}
    \ x_1, x_2, \ldots, x_N\ \mbox{and}\\
    \mbox{ is nondecreasing in}\ t\in [0, \infty) \ \mbox{for every}\ x\in \R^N; \\
    \ h(x, 0)\equiv 0 \ \mbox{for} \ x\in \R^N;\\
      \mbox{there exist constants}\ p \in (2, 2^*) \ \mbox{and}\ \mathcal{C}_1>0 \ \mbox{such that}\\
     |h(x, t)|\le \mathcal{C}_1(1+|t|^{p-2}),\ \ \forall \ (x, t)\in \R^N\times \R.\end{array}\right.\right\}
 \end{equation}

 \vskip4mm
 \par\noindent
 {\bf Corollary 1.9.}\ \ {\it Assume that $V$ and $W$ satisfy }(V0$'$), (V1$'$) {\it and} (W2), and that
 $$
   W(x, u, v)=\sum_{i=1}^{k}\int_{0}^{|\alpha_iu+\beta_iv|}g_i(x, s)s\mathrm{d}s
     +\sum_{j=1}^{l}\int_{0}^{\sqrt{u^2+2\beta_j'uv+\alpha_j'v^2}}h_j(x, s)s\mathrm{d}s,
 $$
 {\it where $\alpha_i, \beta_i, \alpha_j', \beta_j'\in \R$ with $\alpha_i^2+\beta_i^2\ne 0$ and $\alpha_j'>{\beta_j'}^2$,
 $g_i, h_j\in \mathcal{ND}$, $i=1, 2, \ldots, k; j=1, 2, \ldots, l$.
 Then problem \x{HS} has a solution $z^*\in \mathcal{N}^{-}$ such that $\Phi(z^*)=\inf_{\mathcal{N}^{-}}\Phi>0$.}

 \vskip4mm
 \par
   The paper is organized as follows. In the next section, we develop a functional setting to
 deal with \x{HS5} and \x{HS}. Section 3 is devoted to the proof of Theorem 1.8. In Section 4, we discuss the existence of
 ground state solutions of Nehari-Pankov type to \x{HS5}.

 \vskip10mm

 {\section{Variational setting}}
 \setcounter{equation}{0}

 \vskip4mm
 \par
    Let $V_{\varepsilon}(x):=V(\varepsilon x+x_{m})$ and $Q_{\varepsilon}(x):=Q(\varepsilon x+x_{m})$. Then
 we can rewrite \x{HS5} as
 \begin{equation} \label{HS6}
   \left\{
   \begin{array}{ll}
    -\triangle u+ u+V_{\varepsilon}(x)v=Q_{\varepsilon}(x)F_{v}(u, v), \ \ \ \  x\in {\R}^{N},\\
    -\triangle v+ v+V_{\varepsilon}(x)u=Q_{\varepsilon}(x)F_{u}(u, v), \ \ \ \  x\in {\R}^{N},\\
    |u(x)|+|v(x)| \rightarrow 0, \ \ \mbox{as} \ |x|\rightarrow \infty.
   \end{array}
    \right.
 \end{equation}
 We will mainly deal with \x{HS6} instead of \x{HS5}. Let
 $$
   E^{-}=\left\{\left(-\frac{u}{a}, \frac{u}{b}\right) : u\in H^1(\R^N)\right\},
     \ \ \ \   E^{+}=\left\{\left(\frac{u}{a}, \frac{u}{b}\right): u\in H^1(\R^N)\right\}.
 $$
 For any $z=(u, v)\in E$, set
 $$
   z^{-}=\left(\frac{au-bv}{2a}, \frac{bv-au}{2b}\right), \ \ \ \
   z^{+}=\left(\frac{au+bv}{2a}, \frac{au+bv}{2b}\right).
 $$
 It is obvious that $z=z^{-}+z^{+}$. Now we define two new inner products on $E$
 \begin{eqnarray*}
   &     & (z_1, z_2)=\int_{\R^N} \left[\left(\nabla{z}^{+}_1\cdot \nabla{z}^{+}_2 +\nabla{z}^{-}_1\cdot \nabla{z}^{-}_2\right)
              +\left(z_1^{+}\cdot z_2^{+}+z_1^{-}\cdot z_2^{-}\right)\right]\mathrm{d}x,\\
   &     &  \ \ \ \ \ \ \ \ \ \ \ \ \ \ \ \ \ \ \ \ \ \ \ \ \ \ \ \ \ \ \ \ \ \ \ \ \ \ \ \ \ \ \ \ \forall \ z_i=(u_i, v_i)\in E, \ \ i=1, 2
 \end{eqnarray*}
 and
 \begin{eqnarray*}
   &     & (z_1, z_2)_{\dag}=\int_{\R^N} \left[\left(\nabla{z}^{+}_1\cdot \nabla{z}^{+}_2 +\nabla{z}^{-}_1\cdot \nabla{z}^{-}_2\right)
              +V_1(x)\left(z_1^{+}\cdot z_2^{+}+z_1^{-}\cdot z_2^{-}\right)\right]\mathrm{d}x,\\
   &     &  \ \ \ \ \ \ \ \ \ \ \ \ \ \ \ \ \ \ \ \ \ \ \ \ \ \ \ \ \ \ \ \ \ \ \ \ \ \ \ \ \ \ \ \ \forall \ z_i=(u_i, v_i)\in E, \ \ i=1, 2.
 \end{eqnarray*}
 The corresponding norms are
 $$
   \|z\|=\sqrt{(z, z)}, \ \ \ \ \|z\|_{\dag}=\sqrt{(z, z)_{\dag}}, \ \ \ \ \forall \ z=(u, v)\in E.
 $$
 By virtue of (V0) and (V0$'$), it is easy to check that the norms $\|\cdot\|$, $\|\cdot\|_{\dag}$ and  $\|\cdot\|_{H^1(\R^N)}$
 are equivalent on $E$. It is easy to see that $z^{-}$ and $z^{+}$ are orthogonal with respect to the inner products
 $(\cdot, \cdot)$ and $(\cdot, \cdot)_{\dag}$. Thus we have $E=E^{-}\oplus E^{+}$. By a simple calculation, one can get that
 \begin{equation}\label{no}
   \|z\|^2=\int_{\R^N} \left[\left(|\nabla{z}^{+}|^2+|\nabla{z}^{-}|^2\right)+\left(|z^{+}|^2+|z^{-}|^2\right)\right]\mathrm{d}x,
   \ \ \ \ \forall \ z\in E,
 \end{equation}
 \begin{equation}\label{nod}
   \|z\|_{\dag}^2=\int_{\R^N} \left[\left(|\nabla{z}^{+}|^2+|\nabla{z}^{-}|^2\right)+V_1(x)\left(|z^{+}|^2
   +|z^{-}|^2\right)\right]\mathrm{d}x, \ \ \ \ \forall \ z\in E,
 \end{equation}
 $$
   \frac{ab}{a^2+b^2}\left(\|z^{+}\|^2-\|z^{-}\|^2\right)
     =\int_{{\R}^N}\left(\nabla u\cdot \nabla v+uv\right)\mathrm{d}x, \ \ \ \ \forall \ z=(u, v)\in E
 $$
 and
 $$
   \frac{ab}{a^2+b^2}\left(\|z^{+}\|_{\dag}^2-\|z^{-}\|_{\dag}^2\right)
     =\int_{{\R}^N}\left[\nabla u\cdot \nabla v+V_1(x)uv\right]\mathrm{d}x, \ \ \ \ \forall \ z=(u, v)\in E.
 $$
 Therefore, the functionals $\Phi_{\varepsilon}$  defined by \x{Phe} and $\Phi$ by \x{Ph} can be rewritten
 \begin{equation}\label{Phv}
   \Phi_{\varepsilon}(z)=\frac{ab}{a^2+b^2}\left(\|z^{+}\|^2-\|z^{-}\|^2\right)+\frac{1}{2}\int_{\R^{N}}V_{\varepsilon}(x)|z|^2\mathrm{d}x
      -\int_{\R^{N}}Q_{\varepsilon}(x)F(z)\mathrm{d}x, \ \ \ \ \forall \ z\in E
 \end{equation}
 and
 \begin{equation}\label{Ph1}
   \Phi(z)=\frac{ab}{a^2+b^2}\left(\|z^{+}\|_{\dag}^2-\|z^{-}\|_{\dag}^2\right)+\frac{1}{2}\int_{{\R}^N}V_2(x)|z|^2\mathrm{d}x
    -\int_{\R^{N}}W(x,z)\mathrm{d}x, \ \ \ \ \forall \ z\in E,
 \end{equation}
 resectively. Our hypotheses imply that $\Phi_{\varepsilon}, \Phi\in \mathcal{C}^{1}(E, \mathbb{R})$, and a
 standard argument shows that the critical points of $\Phi_{\varepsilon}$ and $\Phi$ are solutions
 of problems \x{PHS} and \x{HS}, respectively. Moreover, by \x{Phde} and \x{Phd}, there hold
 \begin{eqnarray}\label{Phvd1}
   \langle \Phi_{\varepsilon}'(z), \varphi \rangle
    &  =  & \frac{2ab}{a^2+b^2}\left[\left(z^{+}, \varphi^{+}\right)-\left(z^{-}, \varphi^{-}\right)\right]
              +\int_{{\R}^N}V_{\varepsilon}(x)z\cdot \varphi\mathrm{d}x\nonumber\\
    &     &  \ \ \ \    -\int_{{\R}^N}Q_{\varepsilon}(x)F_{z}(z)\cdot \varphi\mathrm{d}x, \ \ \ \ \forall \ z, \ \varphi\in E,
 \end{eqnarray}
 \begin{equation}\label{Phvd2}
   \langle\Phi_{\varepsilon}'(z), z\rangle=\frac{2ab}{a^2+b^2}\left(\|z^{+}\|^2-\|z^{-}\|^2\right)
     +\int_{{\R}^N}V_{\varepsilon}(x)|z|^2\mathrm{d}x
         -\int_{{\R}^N}Q_{\varepsilon}(x)F_{z}(z)\cdot z\mathrm{d}x, \ \ \ \ \forall \ z\in E,
 \end{equation}
 \begin{eqnarray}\label{Phd1}
   \langle \Phi'(z), \varphi \rangle
    &  =  & \frac{2ab}{a^2+b^2}\left[\left(z^{+}, \varphi^{+}\right)_{\dag}-\left(z^{-}, \varphi^{-}\right)_{\dag}\right]
              +\int_{{\R}^N}V_2(x)z\cdot \varphi\mathrm{d}x\nonumber\\
    &     &  \ \ \ \    -\int_{{\R}^N}W_{z}(x, z)\cdot \varphi\mathrm{d}x, \ \ \ \ \forall \ z, \ \varphi\in E
 \end{eqnarray}
 and
 \begin{equation}\label{Phd2}
   \langle\Phi'(z), z\rangle=\frac{2ab}{a^2+b^2}\left(\|z^{+}\|_{\dag}^2-\|z^{-}\|_{\dag}^2\right)+\int_{{\R}^N}V_2(x)|z|^2\mathrm{d}x
         -\langle\Psi'(z), z\rangle, \ \ \ \ \forall \ z\in E.
 \end{equation}

 \vskip4mm
 \par\noindent
 {\bf Lemma 2.1.} {\it Suppose that } (V0), (F1)-(F2) { are satisfied. If $z = (u, v)$ is a critical point of
 $\Phi_{\varepsilon}$, then $|z(x)| \rightarrow 0$ as $|x|\to \infty$. In a word, $z$ is a solution to system
 \x{HS6}.}

 \vskip2mm
 \par
   The proof is almost standard (see \cite[Lemma 2.1 and Theorem 2.1]{dY}.

 \vskip6mm

 {\section{Ground state solutions of Nehari-Pankov type for periodic system}}
 \setcounter{equation}{0}

 \vskip4mm
 \par
   Let $X=X^{-}\oplus X^{+}$ be a real Hilbert space with $X^{-}\bot\  X^{+}$ and $X^{-}$ be separable.
 On $X$ we define a new norm
 \begin{equation}\label{t-n}
    \|u\|_{\tau}:=\max\left\{\|u^{+}\|,\ \sum_{k=1}^{\infty}\frac{1}{2^{k+1}}\left|\left(u^{-}, e_k\right) \right|\right\},
       \ \ \ \ \forall \ u=u^{-}+u^{+}\in X,
 \end{equation}
 where $\{e_k\}_{k=1}^{\infty}$ is a total orthonormal basis of $X^{-}$. The topology generated by
 $\|\cdot\|_{\tau}$ will be denoted by $\tau$ and all topological notions related to it will include the symbol.
 It is clear that
 \begin{equation}\label{t-1}
    \|u^{+}\|\le \|u\|_{\tau}\le \|u\|, \ \ \ \ \forall \ u\in X.
 \end{equation}

   \par
     For a functional $\varphi\in \mathcal{C}^{1}(X, \R)$, $\varphi$ is said to be $\tau$-upper semi-continuous if
 \begin{equation}\label{t-u}
   u_n, u\in X, \ \ \|u_n-u\|_{\tau}\rightarrow 0 \Rightarrow \varphi(u)\ge \limsup_{n\to\infty}\varphi(u_n);
 \end{equation}
 weakly sequentially lower semi-continuous if
 $$
    u_n\rightharpoonup u \ \mbox{in} \ X \Rightarrow \varphi(u)\le \liminf_{n\to\infty}\varphi(u_n);
 $$
 and $\varphi'$ is said to be weakly sequentially continuous if
 $$
   u_n\rightharpoonup u \ \mbox{in} \ X\Rightarrow \lim_{n\to\infty}\langle\varphi'(u_n), v\rangle
     = \langle\varphi'(u), v\rangle, \ \ \ \ \forall \ v\in X.
 $$
 It is easy to see that \x{t-u} holds if and only if
 \begin{equation}\label{t-l}
   u_n, u\in X, \ \ \|u_n-u\|_{\tau}\rightarrow 0 \Rightarrow \varphi(u)\ge \liminf_{n\to\infty}\varphi(u_n).
 \end{equation}

 \vskip4mm
 \par\noindent
 {\bf Lemma 3.1.} (\cite[Theorem 2.4]{Ta1})\ \ {\it Let $X=X^{-}\oplus X^{+}$ be a real Hilbert space with $X^{-}\bot\  X^{+}$
 and $X^{-}$ be separable. Suppose that $\varphi\in \mathcal{C}^{1}(X, \R)$ satisfies the following assumptions:}
 \vskip2mm
 \par
 (H1)\ \ {\it $\varphi$ is $\tau$-upper semi-continuous;}

 \vskip2mm
 \par
 (H2)\ \ {\it $\varphi'$ is weakly sequentially continuous;}

 \vskip2mm
 \par
 (H3)\ \ {\it there exist $r>\rho>0$ and $e\in X^{+}$ with $\|e\|=1$ such that
 $$
   \kappa:=\inf\varphi(S^{+}_{\rho}) > \sup \varphi(\partial \mathfrak{Q}_r),
 $$
 where }
 $$
     S^{+}_{\rho}=\left\{u\in X^{+} : \|u\|=\rho\right\},  \ \ \ \ \mathfrak{Q}_r=\left\{v+se : v\in X^{-},\ s\ge 0,\ \|v+se\|\le r\right\}.
 $$

 \vskip2mm
 \noindent
 {\it Then there exist a constant $c\in [\kappa, \sup \varphi(\mathfrak{Q}_r)]$ and a sequence $\{u_n\}\subset X$ satisfying}
 \begin{equation}\label{Ce0}
   \varphi(u_n)\rightarrow c, \ \ \ \ \|\varphi'(u_n)\|(1+\|u_n\|)\rightarrow 0.
 \end{equation}

 \vskip4mm
 \par
   Let $\Psi(z)=\int_{\R^{N}}W(x,z)\mathrm{d}x$. Employing a standard argument, one can check easily the following lemma.

 \vskip4mm
 \par\noindent
 {\bf Lemma 3.2.}\ \ {\it Suppose that} (V0$'$), (W0) {\it and} (W1) {\it  are satisfied. Then $\Psi$ is nonnegative, weakly sequentially lower
 semi-continuous, and $\Psi'$ is weakly sequentially continuous.}

 \vskip4mm
 \par\noindent
 {\bf Lemma 3.3.}\ \ {\it Suppose that} (V0$'$), (W0), (W1) {\it and} (W4) {\it are satisfied. Then there holds}
 \begin{eqnarray}\label{T301}
   \Phi(z)
    & \ge & \Phi(\theta z+\zeta)+\frac{ab}{a^2+b^2}\|\zeta\|_{\dag}^2-\frac{1}{2}\int_{\R^N}V_2(x)|\zeta|^2\mathrm{d}x\nonumber\\
    &     & \ \ +\frac{1-\theta ^2}{2}\langle\Phi'(z), z \rangle -\theta \langle\Phi'(z), \zeta \rangle,
               \ \ \ \ \forall \ \theta \ge 0, \ \ z\in E, \ \zeta\in E^{-}.
 \end{eqnarray}

 \vskip2mm
 \par\noindent
 {\bf Proof.} \ \ By \x{Ph1}, \x{Phd1}, \x{Phd2} and (W4), one has
 \begin{eqnarray*}
   &     & \Phi(z)-\Phi(\theta z+\zeta)\\
   &  =  & \frac{ab}{a^2+b^2}\|\zeta\|_{\dag}^2-\frac{1}{2}\int_{\R^N}V_2(x)|\zeta|^2\mathrm{d}x
            +\frac{1-\theta^2}{2}\langle\Phi'(z), z \rangle-\theta\langle\Phi'(z), \zeta \rangle\\
   &     & \ \  +\int_{\R^N}\left[\frac{1-\theta^2}{2}\nabla W(x, z)\cdot z-\theta\nabla W(x, z)\cdot \zeta+W(x, \theta z+\zeta)
                -W(x, z)\right]\mathrm{d}x\\
   & \ge & \frac{ab}{a^2+b^2}\|\zeta\|_{\dag}^2-\frac{1}{2}\int_{\R^N}V_2(x)|\zeta|^2\mathrm{d}x
            +\frac{1-\theta^2}{2}\langle\Phi'(z), z \rangle-\theta\langle\Phi'(z), \zeta \rangle,\\
   &     &   \ \ \ \ \ \ \ \ \ \ \ \ \ \ \ \ \ \ \ \ \ \ \ \ \ \ \ \ \ \ \ \ \ \ \ \ \forall \ \theta \ge 0, \ \ z\in E, \ \zeta\in E^{-}.
 \end{eqnarray*}
 This shows that \x{T301} holds.
 \qed

 \vskip4mm
 \par
    From Lemma 3.3, we have the following two corollaries.

 \vskip4mm
 \par\noindent
 {\bf Corollary 3.4.}\ \ {\it Suppose that} (V0$'$), (W0), (W1) {\it and} (W4) {\it are satisfied. Then for
 $z\in \mathcal{N}^{-}$}
 \begin{equation}\label{T401}
   \Phi(z)\ge \Phi(\theta z+\zeta)+\frac{ab}{a^2+b^2}\|\zeta\|_{\dag}^2-\frac{1}{2}\int_{\R^N}V_2(x)|\zeta|^2\mathrm{d}x,
      \ \ \ \ \forall \ \theta\ge 0, \ \ \zeta\in E^{-}.
 \end{equation}

 \vskip4mm
 \par\noindent
 {\bf Corollary 3.5.}\ \ {\it Suppose that} (V0$'$), (W0), (W1) {\it and} (W4) {\it are satisfied. Then}
 \begin{eqnarray}\label{T501}
   \Phi(z)
    & \ge & \frac{ab\theta^2}{a^2+b^2}\|z\|_{\dag}^2 +\frac{\theta^2}{2}\int_{\R^N}V_2(x)\left(|z^{+}|^2-|z^{-}|^2\right)\mathrm{d}x
                 -\int_{{\R}^N}W(x, \theta z^{+})\mathrm{d}x\nonumber\\
    &     & \ \ \ \  +\frac{1-\theta^2}{2}\langle\Phi'(z), z \rangle+\theta^2\langle\Phi'(z), z^{-} \rangle ,\ \ \ \ \forall \ z\in E, \ \ \theta\ge 0.
 \end{eqnarray}

 \vskip4mm
 \par\noindent
 {\bf Lemma 3.6.} \ \ {\it Suppose that} (V0$'$), (W0), (W1) {\it and} (W4) {\it are satisfied. Then}

 \vskip2mm
 \par
  \ \ (i) \ {\it there exists $\rho>0$ such that }
 \begin{eqnarray*}
   m:=\inf_{\mathcal{N}^{-}}\Phi \ge \kappa:=\inf \left\{\Phi(z) : z\in E^{+}, \|z\|_{\dag}=\rho\right\}>0.
 \end{eqnarray*}

 \vskip2mm
 \par
  \ \ (ii) \ {\it $\|z^{+}\|_{\dag}^2\ge \max\left\{\frac{(1-\eta)^2}{2(1+\eta^2)}\|z^{-}\|_{\dag}^2, \frac{(1-\eta)(a^2+b^2)m}{2ab}\right\}$
  for all $z\in \mathcal{N}^{-}$.}

 \vskip2mm
 \par\noindent
 {\bf Proof.} \ \ By (V0$'$), we have $|V_2(x)|\le \frac{2ab\eta}{a^2+b^2}V_1(x)$, it follows from \x{nod} that
 \begin{eqnarray}\label{ab1}
    &     & \frac{2ab}{a^2+b^2}\|z^{+}\|_{\dag}^2 +\int_{\R^N}V_2(x)|z^{+}|^2\mathrm{d}x\nonumber\\
    &  =  & \int_{\R^N} \left[\frac{2ab}{a^2+b^2}|\nabla{z}^{+}|^2
              +\left(\frac{2ab}{a^2+b^2}V_1(x)+V_2(x)\right)|z^{+}|^2\right]\mathrm{d}x\nonumber\\
    & \ge & \frac{2ab}{a^2+b^2}\int_{\R^N} \left[|\nabla{z}^{+}|^2+(1-\eta)V_1(x)|z^{+}|^2\right]\mathrm{d}x\nonumber\\
    & \ge & \frac{2(1-\eta)ab}{a^2+b^2}\|z^{+}\|_{\dag}^2, \ \ \ \ \forall \ z\in E
 \end{eqnarray}
 and
 \begin{eqnarray}\label{ab2}
    &     & \frac{ab}{a^2+b^2}\left(\|z^{+}\|_{\dag}^2-\|z^{-}\|_{\dag}^2\right)+\frac{1}{2}\int_{\R^N} V_2(x)|z|^2\mathrm{d}x\nonumber\\
    &  =  & \frac{ab}{a^2+b^2}\left(\|z^{+}\|_{\dag}^2-\|z^{-}\|_{\dag}^2\right)+\frac{1}{2}\int_{\R^N} V_2(x)\left(|z^{+}|^2
              +|z^{-}|^2\right)\mathrm{d}x+\int_{\R^N} V_2(x)z^{+}\cdot z^{-}\mathrm{d}x\nonumber\\
    & \le & \frac{ab}{a^2+b^2}\left(\|z^{+}\|_{\dag}^2-\|z^{-}\|_{\dag}^2\right)+\frac{1+\eta}{2(1-\eta)}\int_{\R^N} |V_2(x)||z^{+}|^2\mathrm{d}x
              +\frac{1+\eta}{4\eta}\int_{\R^N} |V_2(x)||z^{-}|^2\mathrm{d}x\nonumber\\
    & \le & \frac{ab}{a^2+b^2}\int_{\R^N} \left[|\nabla{z}^{+}|^2+\frac{1+\eta^2}{1-\eta}V_1(x)|z^{+}|^2\right]\mathrm{d}x\nonumber\\
    &     &  \ \  -\frac{ab}{a^2+b^2}\int_{\R^N} \left[|\nabla{z}^{-}|^2+\frac{1-\eta}{2} V_1(x)|z^{-}|^2\right]\mathrm{d}x\nonumber\\
    & \le & \frac{(1+\eta^2)ab}{(1-\eta)(a^2+b^2)}\|z^{+}\|_{\dag}^2-\frac{(1-\eta)ab}{2(a^2+b^2)}\|z^{-}\|_{\dag}^2, \ \ \ \ \forall \ z\in E.
 \end{eqnarray}
 The rest of the proof is standard, so we omit it.
 \qed

 \vskip4mm
 \par\noindent
 {\bf Lemma 3.7.}\ \ {\it Suppose that} (V0$'$), (W0), (W1)  {\it and} (W2) {\it are satisfied.
 Let $e\in E^{+}$ with $\|e\|_{\dag}=1$. Then there is $r_0>\rho$ such that $\sup \Phi(\partial \mathfrak{Q}_r)\le 0$
 for $r\ge r_0$, where}
 \begin{equation}\label{T701}
   \mathfrak{Q}_r=\left\{\zeta+se : \zeta\in E^{-}, s\ge 0, \ \|\zeta+se\|_{\dag}\le r\right\}.
 \end{equation}

 \vskip2mm
 \par\noindent
 {\bf Proof.} \ \ \x{Ph1} and \x{ab2} imply $\Phi(z)\le 0$ for $z\in E^{-}$. Next, it is sufficient to show that
 $\Phi(z)\rightarrow -\infty$ as $z\in E^{-}\oplus \R^{+} e$ and $\|z\|_{\dag}\rightarrow \infty$.
 Arguing indirectly, assume that for some sequence $\{\zeta_n+s_ne\}\subset E^{-}\oplus \R^{+} e$ with
 $\|\zeta_n+s_ne\|_{\dag}\rightarrow \infty$, there is an $M>0$ such that $\Phi(\zeta_n+s_ne)\ge -M$
 for all $n\in \N$. Set $e=\left(\frac{e_0}{a}, \frac{e_0}{b}\right)$, $\zeta_n=\left(-\frac{w_n}{a}, \frac{w_n}{b}\right)$ and $\xi_n=(\zeta_n+s_ne)/\|\zeta_n+s_ne\|_{\dag}
 =\xi_n^{-}+t_ne$, then $\|\xi_n^{-}+t_ne\|_{\dag}=1$. Passing to a subsequence, we may assume that
 $t_n\rightarrow \bar{t}$ and $\xi_n\rightharpoonup \xi$ in $E$, then $\xi_n\rightarrow \xi$ a.e. on $\R^N$,
 $\xi_n^{-}:=\left(-\frac{\tilde{w}_n}{a}, \frac{\tilde{w}_n}{b}\right)\rightharpoonup \xi^{-}:=
 \left(-\frac{\tilde{w}}{a}, \frac{\tilde{w}}{b}\right)$ in $E$. Hence,
 by \x{Ph1} and \x{ab2}, one has
 \begin{eqnarray}\label{T702}
   -\frac{M}{\|\zeta_n+s_ne\|_{\dag}^2}
    & \le & \frac{\Phi(\zeta_n+s_ne)}{\|\zeta_n+s_ne\|_{\dag}^2}\nonumber\\
    &  =  & \frac{abt_n^2}{a^2+b^2}-\frac{ab}{a^2+b^2}\|\xi_n^{-}\|_{\dag}^2+\frac{1}{2}\int_{{\R}^N}V_2(x)|\xi_n^{-}+t_ne|^2\mathrm{d}x\nonumber\\
    &     & \ \ \ \  -\int_{\R^N}\frac{W\left(x, \frac{-w_n+s_ne_0}{a}, \frac{w_n+s_ne_0}{b}\right)}{\|\zeta_n+s_ne\|_{\dag}^2}\mathrm{d}x\nonumber\\
    & \le & \frac{(1+\eta^2)abt_n^2}{(1-\eta)(a^2+b^2)}-\frac{(1-\eta)ab}{2(a^2+b^2)}\|\xi_n^{-}\|_{\dag}^2\nonumber\\
    &     & \ \ \ \  -\int_{\R^N}\frac{W\left(x, \frac{-w_n+s_ne_0}{a}, \frac{w_n+s_ne_0}{b}\right)}{\|\zeta_n+s_ne\|_{\dag}^2}\mathrm{d}x.
 \end{eqnarray}

 \par
   If $\bar{t}=0$, then it follows from \x{T702} that
 \begin{eqnarray*}
   0 & \le & \frac{(1-\eta)ab}{2(a^2+b^2)}\|\xi_n^{-}\|_{\dag}^2
              +\int_{\R^N}\frac{W\left(x, \frac{-w_n+s_ne_0}{a}, \frac{w_n+s_ne_0}{b}\right)}
              {\|\zeta_n+s_ne\|_{\dag}^2}\mathrm{d}x\\
     & \le & \frac{(1+\eta^2)abt_n^2}{(1-\eta)(a^2+b^2)}+\frac{M}{\|\zeta_n+s_ne\|_{\dag}^2}\rightarrow 0,
 \end{eqnarray*}
 which yields $\|\xi_n^{-}\|_{\dag}\rightarrow 0$, and so $1=\|\xi_n\|_{\dag}\rightarrow 0$, a contradiction.

 \par
    If $\bar{t}\ne 0$, then $s_n\rightarrow \infty$. Hence, it follows from \x{T702}, (W2) and Fatou's Lemma that
 \begin{eqnarray*}
  0 & \le & \limsup_{n\to\infty}\left[\frac{(1+\eta^2)abt_n^2}{(1-\eta)(a^2+b^2)}
              -\frac{(1-\eta)ab}{2(a^2+b^2)}\|\xi_n^{-}\|_{\dag}^2\right.\\
    &     &   \ \ \ \ \ \ \ \ \ \ \ \left.-\int_{\R^N}\frac{W\left(x, \frac{-w_n+s_ne_0}{a}, \frac{w_n+s_ne_0}{b}\right)}{\|\zeta_n+s_ne\|_{\dag}^2}\mathrm{d}x\right]\\
    &  =  & \limsup_{n\to\infty}\left[\frac{(1+\eta^2)abt_n^2}{(1-\eta)(a^2+b^2)}
              -\frac{(1-\eta)ab}{2(a^2+b^2)}\|\xi_n^{-}\|_{\dag}^2\right.\\
    &     &   \ \ \ \ \ \ \ \ \ \ \ \left. -\int_{\R^N}\frac{W\left(x, \frac{-w_n+s_ne_0}{a}, \frac{w_n+s_ne_0}{b}\right)}
              {|s_ne_0|^2}|t_ne_0|^2\mathrm{d}x\right]\\
    & \le & \frac{(1+\eta^2)ab}{(1-\eta)(a^2+b^2)}\lim_{n\to\infty}t_n^2-\liminf_{n\to\infty}
               \int_{\R^N}\frac{W\left(x, \frac{-w_n+s_ne_0}{a}, \frac{w_n+s_ne_0}{b}\right)}
              {|s_ne_0|^2}|t_ne_0|^2\mathrm{d}x\\
    & \le & \frac{(1+\eta^2)ab\bar{t}^2}{(1-\eta)(a^2+b^2)}-\int_{\R^N}\liminf_{n\to\infty}
              \left[\frac{W\left(x, \frac{-w_n+s_ne_0}{a}, \frac{w_n+s_ne_0}{b}\right)}{|s_ne_0|^2}|t_ne_0|^2\right]\mathrm{d}x\\
    &  =  & -\infty,
 \end{eqnarray*}
 a contradiction.
 \qed

 \vskip4mm
 \par
    Since $E^{-}$ is separable, let $\{e_k\}_{k=1}^{\infty}$ be a total orthonormal basis of $E^{-}$. On $E$ we
 define the $\tau$-norm
 \begin{equation}\label{tz}
    \|z\|_{\tau}:=\max\left\{\|z^{+}\|_{\dag},\ \sum_{k=1}^{\infty}\frac{1}{2^{k+1}}\left|\left(z^{-}, e_k\right)_{\dag} \right|\right\},
       \ \ \ \ \forall \ z\in E.
 \end{equation}
 It is clear that
 \begin{equation}\label{tz-1}
    \|z^{+}\|_{\dag}\le \|z\|_{\tau}\le \|z\|_{\dag}, \ \ \ \ \forall \ z\in E.
 \end{equation}

 \vskip4mm
 \par\noindent
 {\bf Lemma 3.8.}\ \ {\it Suppose that} (V0$'$), (W0), (W1)  {\it and} (W2)  {\it  are satisfied. Then $\Phi\in \mathcal{C}^{1}(E, \R)$
 is $\tau$-upper semi-continuous and $\Phi'$ is weakly sequentially continuous.}

 \vskip2mm
 \par\noindent
 {\bf Proof.} \ \  It is clear that $\Phi\in \mathcal{C}^{1}(E, \R)$. First we prove that $\Phi$ is $\tau$-upper semi-continuous.
 Let $z_n \xrightarrow {\tau} z$ in $E$ and $\Phi(z_n)\ge c$. It follows from \x{Ph1}, \x{ab2}, \x{tz-1} and (W0) that
 $z_n^{+}\rightarrow z^{+}$ in $E$ and
 \begin{equation*}
   C_1  \ge  \frac{(1+\eta^2)ab}{(1-\eta)(a^2+b^2)}\|z_n^{+}\|_{\dag}^2
                 \ge c+\frac{(1-\eta)ab}{2(a^2+b^2)}\|z_n^{-}\|_{\dag}^2.
 \end{equation*}
 This shows that $\{z_n^{-}\}\subset E^{-}$ is bounded. It is easy to show that $z_n^{-}\xrightarrow{\tau} z^{-}
  \ \Leftrightarrow \  z_n^{-}\rightharpoonup z^{-}$, and so, $z_n\rightarrow z$ a.e. on $\R^N$.
 Note that
 \begin{eqnarray}\label{T704}
    &     & \left|\int_{\R^N} V_2(x)\left(z_n^{+}\cdot z_n^{-}-z^{+}\cdot z^{-}\right)\mathrm{d}x\right|\nonumber\\
    & \le & \int_{\R^N} |V_2(x)|\left|z_n^{+}-z^{+}\right||z_n^{-}|\mathrm{d}x
            +\left|\int_{\R^N} V_2(x)z^{+}\cdot \left(z_n^{-}-z^{-}\right)\mathrm{d}x\right|=o(1)
 \end{eqnarray}
 and
 \begin{eqnarray}\label{T705}
   &     & \liminf_{n\to\infty}\int_{{\R}^N}\left[|\nabla z_n^{-}|^2+\left(V_1(x)
              -\frac{a^2+b^2}{2ab}V_2(x)\right)|z_n^{-}|^2\right]\mathrm{d}x\nonumber\\
   & \ge & \int_{{\R}^N}\left[|\nabla z^{-}|^2+\left(V_1(x)
              -\frac{a^2+b^2}{2ab}V_2(x)\right)|z^{-}|^2\right]\mathrm{d}x.
 \end{eqnarray}
 Hence, it follows from (W2), \x{Ph1}, \x{T704}, \x{T705} and Fatou's Lemma that
 \begin{eqnarray*}
   -\Phi(z)
    &  =  & \frac{ab}{a^2+b^2}\left(\|z^{-}\|_{\dag}^2-\|z^{+}\|_{\dag}^2\right)
              -\frac{1}{2}\int_{{\R}^N}V_2(x)|z|^2\mathrm{d}x+\int_{{\R}^N}W(x, z)\mathrm{d}x\\
    &  =  & \frac{ab}{a^2+b^2}\int_{{\R}^N}\left[|\nabla z^{-}|^2+\left(V_1(x)
              -\frac{a^2+b^2}{2ab}V_2(x)\right)|z^{-}|^2\right]\mathrm{d}x\\
    &     & -\frac{ab}{a^2+b^2}\int_{{\R}^N}\left[|\nabla z^{+}|^2+\left(V_1(x)
              +\frac{a^2+b^2}{2ab}V_2(x)\right)|z^{+}|^2\right]\mathrm{d}x\\
    &     &   -\int_{{\R}^N}V_2(x)z^{+}\cdot z^{-}\mathrm{d}x+\int_{{\R}^N}W(x, z)\mathrm{d}x\\
    & \le & \liminf_{n\to\infty}\left\{\frac{ab}{a^2+b^2}\int_{{\R}^N}\left[|\nabla z_n^{-}|^2+\left(V_1(x)
              -\frac{a^2+b^2}{2ab}V_2(x)\right)|z_n^{-}|^2\right]\mathrm{d}x\right.\\
    &     & -\frac{ab}{a^2+b^2}\int_{{\R}^N}\left[|\nabla z_n^{+}|^2+\left(V_1(x)
              +\frac{a^2+b^2}{2ab}V_2(x)\right)|z_n^{+}|^2\right]\mathrm{d}x\\
    &     &   \left.-\int_{{\R}^N}V_2(x)z_n^{+}\cdot z_n^{-}\mathrm{d}x+\int_{{\R}^N}W(x, z_n)\mathrm{d}x\right\}\\
    &  =  & \liminf_{n\to\infty}[-\Phi(z_n)]= -\limsup_{n\to\infty}\Phi(z_n).
 \end{eqnarray*}
  This shows that $\Phi$ is $\tau$-upper semi-continuous.

 \par
    The proof that $\Phi'$ is weakly sequentially continuous is standard, so we omit it.
 \qed

 \vskip4mm
 \par\noindent
 {\bf Lemma 3.9.}\ \ {\it Suppose that} (V0$'$), (W0), (W1), (W2)  {\it and} (W4) {\it are satisfied.
 Then there exist a constant $c\in [\kappa, \sup \Phi(\mathfrak{Q}_r)]$ for $r\ge r_0\mathfrak{}$ and a sequence $\{z_n\}\subset E$ satisfying
 \begin{equation}\label{T901}
   \Phi(z_n)\rightarrow c, \ \ \ \ \|\Phi'(z_n)\|_{E^*}(1+\|z_n\|_{\dag})\rightarrow 0,
 \end{equation}
 where $\mathfrak{Q}_r$ is defined by \x{T701}.}

 \vskip2mm
 \par\noindent
 {\bf Proof.} \ \ Lemma 3.9 is a direct corollary of Lemmas 3.1, 3.2, 3.6 (i), 3.7 and 3.8.
 \qed

 \vskip4mm
 \par
    Applying Corollary 3.4, Lemmas 3.6 (i), 3.7 and 3.9, we can prove the following lemma in a similar way as
 \cite[Lemma 3.8]{Ta1}.

 \vskip4mm
 \par\noindent
 {\bf Lemma 3.10.}(\cite[Lemma 3.9]{Ta4})\ \ {\it Suppose that} (V0$'$), (W0), (W1), (W2) {\it and} (W4) {\it are satisfied.
 Then there exist a constant $c_*\in [\kappa, m]$ and a sequence $\{z_n\}=\{(u_n, v_n)\}\subset E$ satisfying}
 \begin{equation}\label{Ce}
   \Phi(z_n)\rightarrow c_*, \ \ \ \ \|\Phi'(z_n)\|(1+\|z_n\|)\rightarrow 0.
 \end{equation}

 \vskip4mm
 \par\noindent
 {\bf Lemma 3.11.}\ \ {\it Suppose that} (V0$'$), (W0), (W1), (W2)  {\it and} (W4) {\it are satisfied. Then for any
 $z\in E\setminus E^{-}$, $\mathcal{N}^{-}\cap (E^{-}\oplus \R^{+}z)\ne \emptyset$, i.e.,
 there exist $t(z)>0$ and $\zeta(z)\in E^{-}$ such that $t(z)z+\zeta(z)\in \mathcal{N}^{-}$.}

 \vskip2mm
 \par
    The proof is the same as one of \cite[Lemma 2.6]{SW}.

 \vskip4mm
 \par\noindent
 {\bf Lemma 3.12.}\ \ {\it Suppose that} (V0$'$), (W0), (W1), (W2), (W3) {\it and} (W4) {\it are satisfied. Then any sequence
 $\{z_n\}=\{(u_n, v_n)\}\subset E$ satisfying
 \begin{equation}\label{Ce1}
   \Phi(z_n)\rightarrow c\ge 0, \ \ \ \ \langle\Phi'(z_n), z_n \rangle\rightarrow 0,
      \ \ \ \ \langle\Phi'(z_n), z_n^{-} \rangle\rightarrow 0
 \end{equation}
 is bounded in $E$.}

 \vskip2mm
 \par\noindent
 {\bf Proof.} \ \ To prove the boundedness of $\{z_n\}$, arguing by contradiction, suppose that
 $\|z_n\|_{\dag} \to \infty$. Let $\tilde{z}_n=(\tilde{u}_n, \tilde{v}_n):=z_n/\|z_n\|_{\dag}$. Then $\|\tilde{z}_n\|_{\dag}=1$.
 By Sobolev embedding theorem, there exists a constant $C_2>0$ such that $\|\tilde{z}_n^{+}\|_2 \le C_2$.
 If $\delta:=\limsup_{n\to\infty}\sup_{y\in \R^N}\int_{B_1(y)}|\tilde{z}_n^{+}|^2\mathrm{d}x=0$, then by Lions' concentration compactness
 principle \cite[Lemma 1.21]{Wi}, $\tilde{z}_n^{+}\rightarrow 0$ in $L^{p}(\R^N)$. Fix $\vartheta=[(a^2+b^2)(1+c)/(1-\eta)ab]^{1/2}$.
 By virtue of (W0) and (W1), for $\epsilon =1/4(\vartheta C_2)^2>0$, there exists $C_{\epsilon}>0$ such that
 \begin{eqnarray*}
   W(x, z)\le \epsilon |z|^2+C_{\epsilon}|z|^p,  \ \ \ \ \forall \ (x, z)\in \R^N\times \R^2.
 \end{eqnarray*}
 Hence, it follows that
 \begin{eqnarray}\label{T121}
   \limsup_{n\to\infty}\int_{\R^N}W(x, \vartheta\tilde{z}_n^{+})\mathrm{d}x
   & \le & \epsilon \vartheta^2\limsup_{n\to\infty}\int_{\R^N}|\tilde{z}_n^{+}|^2\mathrm{d}x
             +C_{\epsilon}\vartheta^p\limsup_{n\to\infty}\int_{\R^N}|\tilde{z}_n^{+}|^p\mathrm{d}x\nonumber\\
   & \le & \epsilon(\vartheta C_2)^2=\frac{1}{4}.
 \end{eqnarray}
 Let $\theta_n=\vartheta/\|z_n\|_{\dag}$. Hence, by virtue of \x{T501}, \x{Ce1} and \x{T121}, one can get
 \begin{eqnarray*}
   c+o(1)
   &  =  & \Phi(z_n)\\
   & \ge & \frac{ab\theta_n^2}{a^2+b^2}\|z_n\|_{\dag}^2 +\frac{\theta_n^2}{2}\int_{\R^N}V_2(x)\left(|z_n^{+}|^2-|z_n^{-}|^2\right)\mathrm{d}x
              -\int_{\R^N}W(x, \theta_nz_n^{+})\mathrm{d}x\\
   &     & \ \  +\frac{1-\theta_n^2}{2}\langle\Phi'(z_n), z_n \rangle+\theta_n^2\langle\Phi'(z_n), z_n^{-}\rangle\\
   &  =  & \frac{ab\vartheta^2}{a^2+b^2}\|\tilde{z}_n\|_{\dag}^2 +\frac{\vartheta^2}{2}\int_{\R^N}V_2(x)\left(|\tilde{z}_n^{+}|^2
              -|\tilde{z}_n^{-}|^2\right)\mathrm{d}x -\int_{\R^N}W(x, \vartheta\tilde{z}_n^{+})\mathrm{d}x\\
   &     & \ \  +\left(\frac{1}{2}-\frac{\vartheta^2}{2\|z_n\|_{\dag}^2}\right)\langle\Phi'(z_n), z_n \rangle
                +\frac{\vartheta^2}{\|z_n\|_{\dag}^2}\langle\Phi'(z_n), z_n^{-} \rangle\\
   & \ge & \frac{(1-\eta)ab\vartheta^2}{a^2+b^2} -\int_{\R^N}W(x, \vartheta\tilde{z}_n^{+})\mathrm{d}x +o(1)\\
   & \ge & \frac{(1-\eta)ab\vartheta^2}{a^2+b^2}-\frac{1}{4} +o(1) >  \frac{3}{4}+c+o(1).
 \end{eqnarray*}
 This contradiction shows that $\delta>0$.  Going if necessary to a subsequence, we may assume the existence of $k_n\in \Z^N$ such that
 $\int_{B_{1+\sqrt{N}}(k_n)}|\tilde{z}_n^{+}|^2\mathrm{d}x > \frac{\delta}{2}$. Let $\zeta_n(x):=\tilde{z}_n(x+k_n)$. Then
 \begin{equation}\label{l13}
   \int_{B_{1+\sqrt{N}}(0)}|\zeta_n^{+}|^2\mathrm{d}x> \frac{\delta}{2}.
 \end{equation}
 Now we define $z^{k_n}_n(x):=(u^{k_n}_n(x), v^{k_n}_n(x))=z_n(x+k_n)$, then $z^{k_n}_n/\|z_n\|_{\dag}=\zeta_n$
 and $\|\zeta_n\|_{H^1(\R^N)}^2 =\|\tilde{z}_n\|_{H^1(\R^N)}^2$. Passing to a subsequence, we have $\zeta_n^+\rightharpoonup \zeta^+$
 in $H^1(\R^N)$, $\zeta_n^+\rightarrow \zeta^+$ in $L^{s}_{\mathrm{loc}}(\R^N)$, $2\le s<2^*$ and $\zeta_n^+\rightarrow \zeta^+$ a.e. on $\R^N$.
 Obviously, \x{l13} implies that $\zeta^{+}\ne 0$. For a.e. $x\in \{y\in \R^N : \zeta^{+}(y)\ne 0\}:=\Omega$, we have
 $\lim_{n\to\infty}|au^{k_n}_n(x)+bv^{k_n}_n(x)|=\infty$. Hence, it follows from \x{Ph1}, \x{ab2}, \x{Ce1}, (W2), (W3)
 and Fatou's lemma that
 \begin{eqnarray*}
  0 &  =  & \lim_{n\to\infty}\frac{c+o(1)}{\|z_n\|_{\dag}^2} = \lim_{n\to\infty}\frac{\Phi(z_n)}{\|z_n\|_{\dag}^2}\\
    &  =  & \lim_{n\to\infty}\left[\frac{ab}{a^2+b^2}\left(\|\tilde{z}_n^{+}\|_{\dag}^2-\|\tilde{z}_n^{-}\|_{\dag}^2\right)
             +\frac{1}{2}\int_{\R^N} V_2(x)|\tilde{z}_n|^2\mathrm{d}x
             -\int_{\R^N}\frac{W(x, u_n, v_n)}{\|z_n\|_{\dag}^2}\mathrm{d}x\right]\\
    & \le & \lim_{n\to\infty}\left[\frac{(1+\eta^2)ab}{(1-\eta)(a^2+b^2)}\|\tilde{z}_n^{+}\|_{\dag}^2
             -\frac{(1-\eta)ab}{2(a^2+b^2)}\|\tilde{z}_n^{-}\|_{\dag}^2
             -\int_{\Omega}\frac{W(x+k_n, u^{k_n}_n, v^{k_n}_n)}{\|z_n\|_{\dag}^2}\mathrm{d}x\right]\\
    &  =  & \lim_{n\to\infty}\left[\frac{(1+\eta^2)ab}{(1-\eta)(a^2+b^2)}\|\tilde{z}_n^{+}\|_{\dag}^2
             -\frac{(1-\eta)ab}{2(a^2+b^2)}\|\tilde{z}_n^{-}\|_{\dag}^2
             -\frac{4a^2b^2}{a^2+b^2}\int_{\Omega}\frac{W(x, u^{k_n}_n, v^{k_n}_n)}{|au^{k_n}_n+bv^{k_n}_n|^2}|\zeta_n^{+}|^2\mathrm{d}x\right]\\
    & \le & \frac{(1+\eta^2)ab}{(1-\eta)(a^2+b^2)}-\frac{4a^2b^2}{a^2+b^2}\int_{\Omega}\liminf_{n\to\infty}
             \left[\frac{W(x, u^{k_n}_n, v^{k_n}_n)}{|au^{k_n}_n+bv^{k_n}_n|^2}|\zeta_n^{+}|^2\right]\mathrm{d}x = -\infty.
 \end{eqnarray*}
 This contradiction shows that $\{\|u_n\|_{\dag}\}$ is bounded.
 \qed

 \vskip4mm
 \par
    In the last part of the proof of Lemma 3.12, we make use of the periodicity of $W(x, z)$ on $x$, which is still
 valid by using (W2$'$) instead of (W2) and (W3). Therefore, we have the following lemma.

 \vskip4mm
 \par\noindent
 {\bf Lemma 3.13.}\ \ {\it Suppose that} (V0$'$), (W0), (W1), (W2$'$) {\it and} (W4) {\it are satisfied. Then any sequence
 $\{z_n\}=\{(u_n, v_n)\}\subset E$ satisfying \x{Ce1} is bounded.}

 \vskip4mm
 \par\noindent
 {\bf Lemma 3.14.}(\cite[Lemma 2.3]{Ta3})\ \ {\it Suppose that $t \mapsto h(x, t)$ is nondecreasing on $\R$ and $h(x, 0)=0$ for any $x\in \R^N$.
 Then there holds}
 \begin{equation}\label{L32}
   \left(\frac{1-\theta^2}{2}\tau-\theta\sigma\right)h(x, \tau)|\tau|
      \ge \int_{\theta\tau+\sigma}^{\tau}h(x, s)|s|\mathrm{d}s,  \ \ \ \ \forall \ \theta\ge 0, \ \ \tau,\ \sigma\in \R.
 \end{equation}

  \vskip4mm
 \par\noindent
 {\bf Lemma 3.15.}(\cite{Ta4})\ \ {\it Suppose that $W(x, u, v)=\int_{0}^{|\alpha u+\beta v|}g(x, s)s\mathrm{d}s$,
 where $\alpha, \beta\in \R$ with $\alpha^2+\beta^2\ne 0$ and $g\in \mathcal{ND}$. Then $W$ satisfies} (W0), (W1)
 {\it and} (W4).

 \vskip4mm
 \par\noindent
 {\bf Lemma 3.16.}(\cite{Ta4})\ \ {\it Suppose that $W(x, u, v)=\int_{0}^{\sqrt{u^2+2\beta uv+\alpha v^2}}h(x, s)s\mathrm{d}s$,
 where $\alpha, \beta\in \R$ with $\alpha>\beta^2$ and $h\in \mathcal{ND}$. Then $W$ satisfies} (W0), (W1)
 {\it and} (W4).

 \vskip4mm
 \par\noindent
 {\bf Proof of Theorem 1.8.} \  Applying Lemmas 3.10 and 3.12, we deduce that there exists a bounded sequence
 $\{z_n\}=\{(u_n, v_n)\}\subset E$ satisfying $(\ref{Ce})$. The rest of the proof is standard.
 \qed

 \vskip4mm
 \par
   Employing Theorem 1.8, the conclusion of Corollary 1.9 follows by Lemmas 3.14-3.16.

 \vskip6mm

 {\section{Ground state solutions of Nehari-Pankov type for \x{HS6}}}
 \setcounter{equation}{0}

 \vskip4mm
 \par
    Without loss of generality, from now on, we assume that $x_{v}=0\in \mathcal{V}$. We only consider the case when
 (V1) is satisfied, since the arguments are similar when (V2) is satisfied. Then
 \begin{equation}\label{V0}
   V(0) = V_{\min}, \ \ \ \ Q(x)\le Q(0), \ \ \ \ \forall \ |x|\ge R.
 \end{equation}
    Let $V_1=1$ and $V_2=V_{\varepsilon}$ (or $\hat{V}, V_{\min}, V_{\max}$), $W(x, z)=Q_{\varepsilon}(x)F(z)$
 (or $Q(0)F(z), Q_{\min}F(z)$). Then (V0), (F1) and (F2) imply (V0$'$), (W0), (W1), (W2$'$)  and (W4), respectively.
 Let
 $$
    \hat{V}:=\frac{1}{2}(V_{\infty}+V_{\min})=\frac{1}{2}(V_{\infty}+V(0)).
 $$

 \par
    We define three auxiliary functionals as follows:
 \begin{equation}\label{Phh}
   \hat{\Phi}(z)=\frac{ab}{a^2+b^2}\left(\|z^{+}\|^2-\|z^{-}\|^2\right)+\frac{\hat{V}}{2}\int_{\R^{N}}|z|^2\mathrm{d}x
      -Q(0)\int_{\R^{N}}F(z)\mathrm{d}x, \ \ \ \ \forall \ z\in E,
 \end{equation}
 \begin{equation}\label{Ph0}
   \Phi_0(z)=\frac{ab}{a^2+b^2}\left(\|z^{+}\|^2-\|z^{-}\|^2\right)+\frac{V(0)}{2}\int_{\R^{N}}|z|^2\mathrm{d}x
      -Q(0)\int_{\R^{N}}F(z)\mathrm{d}x, \ \ \ \ \forall \ z\in E
 \end{equation}
 and
 \begin{equation}\label{Ph*}
   \Phi_{*}(z)=\frac{ab}{a^2+b^2}\left(\|z^{+}\|^2-\|z^{-}\|^2\right)+\frac{V_{\max}}{2}\int_{\R^{N}}|z|^2\mathrm{d}x
      -Q_{\min}\int_{\R^{N}}F(z)\mathrm{d}x, \ \ \ \ \forall \ z\in E.
 \end{equation}
 Let
 \begin{equation}\label{Neh-}
   \mathcal{\hat{N}}^{-}  = \left\{z\in E\setminus E^{-} : \langle \hat{\Phi}'(z), z \rangle=\langle \hat{\Phi}'(z), \zeta \rangle=0,
       \ \forall \ \zeta\in E^{-} \right\}
 \end{equation}
 and
 \begin{equation}\label{Ne0-}
   \mathcal{N}_0^{-}  = \left\{z\in E\setminus E^{-} : \langle \Phi_0'(z), z \rangle=\langle \Phi_0'(z), \zeta \rangle=0,
       \ \forall \ \zeta\in E^{-} \right\}
 \end{equation}
 be the Nehari-Pankov ``manifolds" of the functionals $\hat{\Phi}$ and $\Phi_0$, respectively. Let
 \begin{equation}\label{cc}
   c_{\varepsilon}=\inf_{\mathcal{N}_{\varepsilon}^{-}}\Phi_{\varepsilon}, \ \ \ \
   \hat{c}=\inf_{\mathcal{\hat{N}}^{-}}\hat{\Phi}, \ \ \ \
   c_{0}=\inf_{\mathcal{N}_{0}^{-}}\Phi_{0}.
 \end{equation}

 \vskip4mm
 \par
   Applying Lemma 3.3 and Corollary 3.5 to $\Phi_{\varepsilon}$, we have the following two lemmas.

 \vskip4mm
 \par\noindent
 {\bf Lemma 4.1.}\ \ {\it Suppose that} (V0) {\it and} (F1) {\it  are satisfied. Then }
 \begin{eqnarray}\label{L01}
   \Phi_{\varepsilon}(z)
    & \ge & \Phi_{\varepsilon}(\theta z+\zeta)+\frac{ab}{a^2+b^2}\|\zeta\|^2-\frac{1}{2}\int_{\R^{N}}V_{\varepsilon}(x)|\zeta|^2\mathrm{d}x
             +\frac{1-\theta ^2}{2}\langle\Phi_{\varepsilon}'(z), z \rangle\nonumber\\
    &     &   -\theta \langle\Phi_{\varepsilon}'(z), \zeta \rangle, \ \ \ \ \forall \ \theta \ge 0, \ \ z\in E, \ \ \zeta\in E^{-}\\
    & \ge & \Phi_{\varepsilon}(\theta z+\zeta)+\frac{(1-\eta)ab}{a^2+b^2}\|\zeta\|^2
             +\frac{1-\theta ^2}{2}\langle\Phi_{\varepsilon}'(z), z \rangle
             -\theta \langle\Phi_{\varepsilon}'(z), \zeta \rangle, \nonumber\\
    &     & \ \ \ \ \forall \ \theta \ge 0, \ \ z\in E, \ \ \zeta\in E^{-}.  \label{L02}
 \end{eqnarray}

 \vskip4mm
 \par\noindent
 {\bf Lemma 4.2.}\ \ {\it Suppose that} (V0) {\it and} (F1) {\it  are satisfied. Then }
 \begin{eqnarray}\label{L21}
   \Phi_{\varepsilon}(z)
    & \ge & \frac{ab\theta^2}{a^2+b^2}\|z\|^2+\frac{\theta^2}{2}\int_{\R^{N}}V_{\varepsilon}(x)\left(|z^{+}|^2-|z^{-}|^2\right)\mathrm{d}x
             -\int_{\R^{N}}Q_{\varepsilon}(x)F(\theta z^{+})\mathrm{d}x\nonumber\\
    &     & +\frac{1-\theta^2}{2}\langle\Phi_{\varepsilon}'(z), z \rangle +\theta^2 \langle\Phi_{\varepsilon}'(z), z^{-} \rangle,
              \ \ \ \ \forall \ \theta \ge 0, \ \ z\in E\\
    & \ge & \frac{(1-\eta)ab\theta^2}{a^2+b^2}\|z\|^2-\int_{\R^{N}}Q_{\varepsilon}(x)F(\theta z^{+})\mathrm{d}x\nonumber\\
    &     &  +\frac{1-\theta ^2}{2}\langle\Phi_{\varepsilon}'(z), z \rangle
             +\theta^2 \langle\Phi_{\varepsilon}'(z), z^{-} \rangle, \ \ \ \ \forall \ \theta \ge 0, \ \ z\in E. \label{L22}
 \end{eqnarray}

 \vskip4mm
 \par
   By virtue of Corollary 1.9, under assumptions (V0), (F1) and (F2), there exists a $\hat{z}\in \mathcal{\hat{N}}^{-}$ such that
 $\hat{c}=\hat{\Phi}(\hat{z})$. In view of Lemma 3.11, there exist $\hat{t}> 0$ and $\hat{\zeta}\in E^{-}$ such that
 $\hat{t}\hat{z}+\hat{\zeta}\in \mathcal{N}_0^{-}$, and so $\Phi_0(\hat{t}\hat{z}+\hat{\zeta})\ge c_0$.

 \vskip4mm
 \par\noindent
 {\bf Lemma 4.3.}\ \ {\it Suppose that} (V0), (V1), (F1) {\it and} (F2)  {\it are satisfied. Then
 $ \hat{c} \ge c_0+\hat{\delta}$, where
 \begin{equation}\label{L34}
   \hat{\delta}:=\frac{V_{\infty}-V_{\min}}{4}\int_{\R^{N}}|\hat{t}\hat{z}+\hat{\zeta}|^2\mathrm{d}x>0
 \end{equation}
 is independent of $\varepsilon>0$.}

 \vskip2mm
 \par\noindent
 {\bf Proof.} \ \ Applying Lemma 3.3 to $\hat{\Phi}(z)$, one has
 \begin{eqnarray*}\label{L33}
   \hat{c}
    &  =  & \hat{\Phi}(\hat{z})\ge \hat{\Phi}(\hat{t}\hat{z}+\hat{\zeta})\nonumber\\
    &  =  & \Phi_0(\hat{t}\hat{z}+\hat{\zeta})+\frac{\hat{V}-V_{\min}}{2}\int_{\R^{N}}|\hat{t}\hat{z}+\hat{\zeta}|^2\mathrm{d}x\nonumber\\
    & \ge & c_0+\frac{V_{\infty}-V_{\min}}{4}\int_{\R^{N}}|\hat{t}\hat{z}+\hat{\zeta}|^2\mathrm{d}x = c_0+\hat{\delta}.
 \end{eqnarray*}
 \qed

 \par
    By virtue of Corollary 1.9, under assumptions (V0), (F1) and (F2), there exists a $z_0\in \mathcal{N}_0^{-}$
 such that $c_0=\Phi_0(z_0)$. Then
 \begin{equation}\label{L40}
   \Phi_0(z_0)\ge \Phi_0(tz_0+\zeta), \ \ \ \ \langle\Phi_0'(z_0), z_0\rangle=\langle\Phi_0'(z_0), \zeta\rangle=0,
     \ \ \ \ \forall \ t\ge 0, \ \ \zeta\in E^{-}.
 \end{equation}
 In view of Lemma 3.11, for any $\varepsilon>0$, there exist $t_{\varepsilon}> 0$ and $\zeta_{\varepsilon}\in E^{-}$
 such that $t_{\varepsilon}z_0+\zeta_{\varepsilon}\in \mathcal{N}_{\varepsilon}^{-}$, and so
 $\Phi_{\varepsilon}(t_{\varepsilon}z_0+\zeta_{\varepsilon})\ge c_{\varepsilon}$ and
 $\Phi_{\varepsilon}(t_{\varepsilon}z_0+\zeta_{\varepsilon})\ge \Phi_{\varepsilon}(tz_0+\zeta), \ \forall \ t\ge 0,
 \ \zeta\in E^{-}$. Set
 \begin{equation}\label{al0}
    \alpha_0:=\frac{ab}{a^2+b^2}\left(\|z_0^{+}\|^2-\|z_0^{-}\|^2\right)+\frac{V_{\min}}{2}\int_{\R^{N}}|z_0|^2\mathrm{d}x
             -Q_{\max}\int_{\R^{N}}F(z_0)\mathrm{d}x.
 \end{equation}
 Clearly, $\alpha_0$ is independent of $\varepsilon>0$. Analogous to the proof of Lemma 3.7, one can demonstrate
 the following lemma.

 \vskip4mm
 \par\noindent
 {\bf Lemma 4.4.}\ \ {\it Suppose that} (V0), (V1), (F1) {\it and} (F2)  {\it are satisfied.
 Then there is an $M_0>0$ independent of $\varepsilon>0$ such that}
 \begin{equation}\label{r0}
   \Phi_*(\zeta+sz_0)\le \alpha_0-1, \ \ \ \ \forall \ \zeta\in E^{-}, \ s\ge 0, \ \|\zeta+sz_0\|\ge M_0.
 \end{equation}

 \vskip4mm
 \par\noindent
 {\bf Lemma 4.5.}\ \ {\it Suppose that} (V0), (V1), (F1) {\it and} (F2)  {\it are satisfied. Then}
 \begin{equation}\label{M12}
   M_1:=\sup_{\varepsilon>0}|t_{\varepsilon}| \le M_0\|z_0^{+}\|^{-1}, \ \ \ \ M_2:=\sup_{\varepsilon>0}\|\zeta_{\varepsilon}\|\le M_0(1+\|z_0\|\|z_0^{+}\|^{-1}).
 \end{equation}

 \vskip2mm
 \par\noindent
 {\bf Proof.} \ \ Note that
 \begin{eqnarray}\label{L44}
   \alpha_0
    &  =  & \frac{ab}{a^2+b^2}\left(\|z_0^{+}\|^2-\|z_0^{-}\|^2\right)+\frac{V_{\min}}{2}\int_{\R^{N}}|z_0|^2\mathrm{d}x
             -Q_{\max}\int_{\R^{N}}F(z_0)\mathrm{d}x\nonumber\\
    & \le & \Phi_{\varepsilon}(z_0)\le\Phi_{\varepsilon}(t_{\varepsilon}z_0+\zeta_{\varepsilon})
             \le \Phi_{*}(t_{\varepsilon}z_0+\zeta_{\varepsilon}).
 \end{eqnarray}
 By Lemma 4.4, one obtain that
 \begin{equation}\label{L42}
   \sup_{\varepsilon > 0}\|t_{\varepsilon}z_0+\zeta_{\varepsilon}\|\le M_0.
 \end{equation}
 Since
 $$
   \|t_{\varepsilon}z_0+\zeta_{\varepsilon}\|^2=\|t_{\varepsilon}z_0^{+}\|^2+\|t_{\varepsilon}z_0^{-}+\zeta_{\varepsilon}\|^2\ge t_{\varepsilon}^2\|z_0^{+}\|^2,
 $$
 it follows from \x{L42} and the above that
 \begin{equation*}
   M_1=\sup_{\varepsilon>0}|t_{\varepsilon}| \le M_0\|z_0^{+}\|^{-1}, \ \ \ \ M_2=\sup_{\varepsilon>0}\|\zeta_{\varepsilon}\|\le M_0(1+\|z_0\|\|z_0^{+}\|^{-1}).
 \end{equation*}
 \qed

 \par
    In view of (F1), there exists a constant $\beta_0>0$ such that
 \begin{equation}\label{Fz}
   |F_{z}(z)| \le \beta_0(|z|+|z|^{p-1}), \ \ \ \ \forall \ z\in \R^2.
 \end{equation}

 \par
    Now, we can choose $R_0>R$ sufficient large such that
 \begin{equation}\label{VQR}
   V(x)\ge \hat{V}, \ \ \ \ Q(x)\le Q(0), \ \ \ \ \forall \ |x|\ge R_0\
 \end{equation}
 and
 \begin{eqnarray}\label{L58}
    &     & M_1^2\int_{|x|> R_0}|z_0|^2\mathrm{d}x-(1+M_1^2)Q_{\max}\int_{|x|> R_0}F_{z}(z_0)\cdot z_0\mathrm{d}x\nonumber\\
    &     &  +2(1+\beta_0Q_{\max})M_1M_2\left(\int_{|x|> R_0}|z_0|^2\mathrm{d}x\right)^{1/2}\nonumber\\
    &     &  +2\beta_0\gamma_pQ_{\max}M_1M_2\left(\int_{|x|> R_0}|z_0|^p\mathrm{d}x\right)^{(p-1)/p} \le \frac{\hat{\delta}}{2},
 \end{eqnarray}
 where $\gamma_p$ is the embedding constant with $\|\cdot\|_{p}\le \gamma_p\|\cdot\|$. For the $R_0>0$ given above, we can
 choose an $\varepsilon_0>0$ such that
 \begin{eqnarray}\label{L60}
    &     & \sqrt{2+2\beta_0Q_{\max}}M_1M_2\|z_0\|_2\left\{\sup_{|x|\le R_0}\left[\left|V_{\varepsilon}(x)-V(0)\right|
             +\beta_0\left|Q_{\varepsilon}(x)-Q(0)\right|\right]\right\}^{1/2}\nonumber\\
    &     & +\frac{M_1^2}{2}\sup_{|x|\le R_0}\left|V_{\varepsilon}(x)-V(0)\right|\|z_0\|_2^2
              +\frac{1+M_1^2}{2}\sup_{|x|\le R_0}\left|Q_{\varepsilon}(x)-Q(0)\right|\int_{|x|\le R_0}F_{z}(z_0)\cdot z_0\mathrm{d}x\nonumber\\
    &     & +\beta_0\gamma_p(2Q_{\max})^{1/p}M_1M_2\|z_0\|_p^{p-1}\left\{\sup_{|x|\le R_0}\left|Q_{\varepsilon}(x)-Q(0)\right|\right\}^{(p-1)/p}
              \le \frac{\hat{\delta}}{4}, \ \ \ \ \varepsilon\in [0, \varepsilon_0]. \ \ \ \ \ \
 \end{eqnarray}

 \vskip4mm
 \par\noindent
 {\bf Lemma 4.6.}\ \ {\it Suppose that} (V0), (V1), (F1) {\it and} (F2)  {\it  are satisfied. Then}
 \begin{equation}\label{c0v}
   c_{0} \ge c_{\varepsilon}-3\hat{\delta}/4, \ \ \ \ \forall \ \varepsilon\in [0, \varepsilon_0].
 \end{equation}

 \vskip2mm
 \par\noindent
 {\bf Proof.} \ \ From (F1), \x{Phv}, \x{Ph0}, \x{L02}, \x{L40}, \x{M12}, \x{Fz}, \x{L58}, \x{L60} and the H\"older inequality, we have
 \begin{eqnarray*}\label{L61}
   c_0
    &  =  & \Phi_0(z_0)=\Phi_{\varepsilon}(z_0)+\frac{1}{2}\int_{\R^{N}}\left[V(0)-V_{\varepsilon}(x)\right]|z_0|^2\mathrm{d}x
              +\int_{\R^{N}}\left[Q_{\varepsilon}(x)-Q(0)\right]F(z_0)\mathrm{d}x\nonumber\\
    & \ge & \Phi_{\varepsilon}(t_{\varepsilon}z_0+\zeta_{\varepsilon})+\frac{(1-\eta)ab}{a^2+b^2}\|\zeta_{\varepsilon}\|^2
             +\frac{1-t_{\varepsilon}^2}{2}\langle\Phi_{\varepsilon}'(z_0), z_0\rangle
             -t_{\varepsilon}\langle\Phi_{\varepsilon}'(z_0), \zeta_{\varepsilon}\rangle\nonumber\\
    &     & +\frac{1}{2}\int_{\R^{N}}\left[V(0)-V_{\varepsilon}(x)\right]|z_0|^2\mathrm{d}x
              +\int_{\R^{N}}\left[Q_{\varepsilon}(x)-Q(0)\right]F(z_0)\mathrm{d}x\nonumber\\
    & \ge & c_{\varepsilon}+\frac{1}{2}\int_{\R^{N}}\left[V(0)-V_{\varepsilon}(x)\right]|z_0|^2\mathrm{d}x
              +\int_{\R^{N}}\left[Q_{\varepsilon}(x)-Q(0)\right]F(z_0)\mathrm{d}x\nonumber\\
    &     & +\frac{1-t_{\varepsilon}^2}{2}\left\{\int_{\R^{N}}\left[V_{\varepsilon}(x)-V(0)\right]|z_0|^2\mathrm{d}x
             +\int_{\R^{N}}\left[Q(0)-Q_{\varepsilon}(x)\right]F_{z}(z_0)\cdot z_0\mathrm{d}x \right\}
             \nonumber\\
    &     & -t_{\varepsilon}\left\{\int_{\R^{N}}\left[V_{\varepsilon}(x)-V(0)\right]z_0\cdot\zeta_{\varepsilon}\mathrm{d}x
             +\int_{\R^{N}}\left[Q(0)-Q_{\varepsilon}(x)\right]F_{z}(z_0)\cdot \zeta_{\varepsilon}\mathrm{d}x\right\}\nonumber\\
    &  =  & c_{\varepsilon}+\frac{1}{2}\int_{\R^{N}}\left[Q_{\varepsilon}(x)-Q(0)\right][2F(z_0)-F_{z}(z_0)\cdot z_0]\mathrm{d}x\nonumber\\
    &     &  -\frac{t_{\varepsilon}^2}{2}\left\{\int_{\R^{N}}\left[V_{\varepsilon}(x)-V(0)\right]|z_0|^2\mathrm{d}x
             +\int_{\R^{N}}\left[Q(0)-Q_{\varepsilon}(x)\right]F_{z}(z_0)\cdot z_0\mathrm{d}x \right\}\nonumber\\
    &     & -t_{\varepsilon}\left\{\int_{\R^{N}}\left[V_{\varepsilon}(x)-V(0)\right]z_0\cdot\zeta_{\varepsilon}\mathrm{d}x
             +\int_{\R^{N}}\left[Q(0)-Q_{\varepsilon}(x)\right]F_{z}(z_0)\cdot \zeta_{\varepsilon}\mathrm{d}x\right\}\nonumber\\
    & \ge & c_{\varepsilon}-\frac{1}{2}\int_{\R^{N}}\left|Q_{\varepsilon}(x)-Q(0)\right||F_{z}(z_0)\cdot z_0-2F(z_0)|\mathrm{d}x\nonumber\\
    &     &  -\frac{M_1^2}{2}\int_{\R^{N}}\left|V_{\varepsilon}(x)-V(0)\right||z_0|^2\mathrm{d}x
             -\frac{M_1^2}{2}\int_{\R^{N}}\left|Q_{\varepsilon}(x)-Q(0)\right|F_{z}(z_0)\cdot z_0\mathrm{d}x\nonumber\\
    &     &  -M_1\int_{\R^{N}}\left|V_{\varepsilon}(x)-V(0)\right||z_0||\zeta_{\varepsilon}|\mathrm{d}x
             -M_1\int_{\R^{N}}\left|Q_{\varepsilon}(x)-Q(0)\right||F_{z}(z_0)||\zeta_{\varepsilon}|\mathrm{d}x\nonumber\\
    & \ge & c_{\varepsilon}-\frac{M_1^2}{2}\int_{\R^{N}}\left|V_{\varepsilon}(x)-V(0)\right||z_0|^2\mathrm{d}x
             -\frac{1+M_1^2}{2}\int_{\R^{N}}\left|Q_{\varepsilon}(x)-Q(0)\right|F_{z}(z_0)\cdot z_0\mathrm{d}x\nonumber\\
    &     & -\sqrt{2+2\beta_0Q_{\max}}M_1M_2\left\{\int_{\R^{N}}\left[\left|V_{\varepsilon}(x)-V(0)\right|
             +\beta_0\left|Q_{\varepsilon}(x)-Q(0)\right|\right]|z_0|^2\mathrm{d}x\right\}^{1/2}\nonumber\\
    &     &  -\beta_0\gamma_p(2Q_{\max})^{1/p}M_1M_2\left\{\int_{\R^{N}}\left|Q_{\varepsilon}(x)
              -Q(0)\right||z_0|^p\mathrm{d}x\right\}^{(p-1)/p}\nonumber\\
    & \ge & c_{\varepsilon}-\frac{M_1^2}{2}\int_{|x|\le R_0}\left|V_{\varepsilon}(x)-V(0)\right||z_0|^2\mathrm{d}x
              -\frac{1+M_1^2}{2}\int_{|x|\le R_0}\left|Q_{\varepsilon}(x)-Q(0)\right|F_{z}(z_0)\cdot z_0\mathrm{d}x\nonumber\\
    &     & -\sqrt{2+2\beta_0Q_{\max}}M_1M_2\left\{\int_{|x|\le R_0}\left[\left|V_{\varepsilon}(x)-V(0)\right|
             +\beta_0\left|Q_{\varepsilon}(x)-Q(0)\right|\right]|z_0|^2\mathrm{d}x\right\}^{1/2}\nonumber\\
    &     & -\beta_0\gamma_p(2Q_{\max})^{1/p}M_1M_2\left\{\int_{|x|\le R_0}\left|Q_{\varepsilon}(x)
               -Q(0)\right||z_0|^p\mathrm{d}x\right\}^{(p-1)/p}\nonumber\\
    &     & -M_1^2\int_{|x|> R_0}|z_0|^2\mathrm{d}x-(1+M_1^2)Q_{\max}\int_{|x|> R_0}F_{z}(z_0)\cdot z_0\mathrm{d}x\nonumber\\
    &     &  -2(1+\beta_0Q_{\max})M_1M_2\left(\int_{|x|> R_0}|z_0|^2\mathrm{d}x\right)^{1/2}
             -2\beta_0\gamma_pQ_{\max}M_1M_2\left(\int_{|x|> R_0}|z_0|^p\mathrm{d}x\right)^{(p-1)/p}\nonumber\\
    & \ge & c_{\varepsilon}-\sqrt{2+2\beta_0Q_{\max}}M_1M_2\|z_0\|_2\left\{\sup_{|x|\le R_0}\left[\left|V_{\varepsilon}(x)-V(0)\right|
             +\beta_0\left|Q_{\varepsilon}(x)-Q(0)\right|\right]\right\}^{1/2}\nonumber\\
    &     & -\frac{M_1^2}{2}\sup_{|x|\le R_0}\left|V_{\varepsilon}(x)-V(0)\right|\|z_0\|_2^2
              -\frac{1+M_1^2}{2}\sup_{|x|\le R_0}\left|Q_{\varepsilon}(x)-Q(0)\right|\int_{|x|\le R_0}F_{z}(z_0)\cdot z_0\mathrm{d}x\nonumber\\
    &     & -\beta_0\gamma_p(2Q_{\max})^{1/p}M_1M_2\|z_0\|_p^{p-1}\left\{\sup_{|x|\le R_0}\left|Q_{\varepsilon}(x)-Q(0)\right|\right\}^{(p-1)/p}
            -\frac{\hat{\delta}}{2}\nonumber\\
    & \ge & c_{\varepsilon}-\frac{3\hat{\delta}}{4}.
 \end{eqnarray*}
 \qed

 \vskip4mm
 \par
   Similar to Lemma 3.6, we can demonstrate that for any $\varepsilon>0$, there exists a $\rho_{\varepsilon}>0$ such that
 \begin{eqnarray}\label{Ka}
   c_{\varepsilon}=\inf_{\mathcal{N}_{\varepsilon}^{-}}\Phi_{\varepsilon} \ge \kappa_{\varepsilon}:
     =\inf \left\{\Phi_{\varepsilon}(z) : z\in E^{+}, \|z\|=\rho_{\varepsilon}\right\}>0.
 \end{eqnarray}
 Applying Lemmas 3.10 and 3.13 to $\Phi_{\varepsilon}$, we have the following two lemmas.

 \vskip4mm
 \par\noindent
 {\bf Lemma 4.7.}\ \ {\it Suppose that} (V0), (F1) {\it and} (F2) {\it  are satisfied. Then
 there exist a constant $\bar{c}_{\varepsilon}\in [\kappa_{\varepsilon}, c_{\varepsilon}]$ and a sequence $\{z_n^{\varepsilon}\}
 =\{(u_n^{\varepsilon}, v_n^{\varepsilon})\}\subset E$ satisfying}
 \begin{equation}\label{Ce3}
   \Phi_{\varepsilon}(z_n^{\varepsilon}) \rightarrow \bar{c}_{\varepsilon},
    \ \ \ \ \|\Phi_{\varepsilon}'(z_n^{\varepsilon})\|(1+\|z_n^{\varepsilon}\|)\rightarrow 0.
 \end{equation}

 \vskip4mm
 \par\noindent
 {\bf Lemma 4.8.}\ \ {\it Suppose that} (V0), (F1) {\it and} (F2) {\it are satisfied. Then the sequence
 $\{z_n^{\varepsilon}\}=\{(u_n^{\varepsilon}, v_n^{\varepsilon})\}\subset E$ satisfying \x{Ce3} is bounded in $E$.}

 \vskip4mm
 \par
   Similar to \x{ab2}, one has
 \begin{eqnarray}\label{ab3}
    &     & \frac{ab}{a^2+b^2}\left(\|z^{+}\|^2-\|z^{-}\|^2\right)+\frac{1}{2}\int_{\R^N} V_{\varepsilon}(x)|z|^2\mathrm{d}x\nonumber\\
    & \le & \frac{(1+\eta^2)ab}{(1-\eta)(a^2+b^2)}\|z^{+}\|^2-\frac{(1-\eta)ab}{2(a^2+b^2)}\|z^{-}\|^2, \ \ \ \ \forall \ z\in E.
 \end{eqnarray}
 By (F1), \x{Phv}, \x{Ce3} and \x{ab3}, one has
 \begin{eqnarray}\label{ab4}
   \bar{c}_{\varepsilon}+o(1)
    &  =  & \Phi_{\varepsilon}(z_n^{\varepsilon})\nonumber\\
    &  =  & \frac{ab}{a^2+b^2}\left(\|{z_n^{\varepsilon}}^{+}\|^2-\|{z_n^{\varepsilon}}^{-}\|^2\right)
             +\frac{1}{2}\int_{\R^{N}}V_{\varepsilon}(x)|z_n^{\varepsilon}|^2\mathrm{d}x
             -\int_{\R^{N}}Q_{\varepsilon}(x)F(z_n^{\varepsilon})\mathrm{d}x\nonumber\\
    & \le & \frac{(1+\eta^2)ab}{(1-\eta)(a^2+b^2)}\|{z_n^{\varepsilon}}^{+}\|^2
             -\frac{(1-\eta)ab}{2(a^2+b^2)}\|{z_n^{\varepsilon}}^{-}\|^2\nonumber\\
    & \le & \frac{(1+\eta^2)ab}{(1-\eta)(a^2+b^2)}\|{z_n^{\varepsilon}}^{+}\|^2.
 \end{eqnarray}
 Therefore, it follows from Lemma 3.11 that there exist $t_n^{\varepsilon}>0$ and $\zeta_n^{\varepsilon}\in E^{-}$ such that
 $t_n^{\varepsilon}z_n^{\varepsilon}+\zeta_n^{\varepsilon}\in \hat{\mathcal{N}}^{-}$, and so
 \begin{equation}\label{L85}
   \hat{\Phi}(t_n^{\varepsilon}z_n^{\varepsilon}+\zeta_n^{\varepsilon})\ge \hat{c}, \ \ \ \
   \langle \hat{\Phi}'(t_n^{\varepsilon}z_n^{\varepsilon}+\zeta_n^{\varepsilon}), t_n^{\varepsilon}z_n^{\varepsilon}+\zeta_n^{\varepsilon} \rangle
    =\langle \hat{\Phi}'(t_n^{\varepsilon}z_n^{\varepsilon}+\zeta_n^{\varepsilon}), \zeta \rangle=0, \ \ \ \ \forall \ \zeta\in E^{-}.
 \end{equation}

 \vskip4mm
 \par\noindent
 {\bf Lemma 4.9.}\ \ {\it Suppose that} (V0), (V1) {\it and} (F1)-(F3) {\it  are satisfied. Then for any $\varepsilon>0$, there exist
 $K_1(\varepsilon)>0$ and $K_2(\varepsilon)>0$ such that}
 \begin{equation}\label{tnva}
   0\le t_n^{\varepsilon}\le K_1(\varepsilon), \ \ \ \ \|\zeta_n^{\varepsilon}\|\le K_2(\varepsilon), \ \ \ \ \forall \ \varepsilon> 0.
 \end{equation}

 \vskip2mm
 \par\noindent
 {\bf Proof.} \ \ If along a subsequence $t_n^{\varepsilon}<T_0$, we are through. So we may assume that
 $t_n^{\varepsilon}\ge T_0$. In view of Lemma 4.8, there exists a constant $C_3>0$ such that $\|{z_n^{\varepsilon}}^{+}\|_2\le C_3$.
 If $\delta_{\varepsilon}:=\limsup_{n\to\infty}\sup_{y\in \R^N}\int_{B_1(y)}|{z_n^{\varepsilon}}^{+}|^2\mathrm{d}x=0$, then by
 Lions' concentration compactness principle \cite[Lemma 1.21]{Wi}, ${z_n^{\varepsilon}}^{+}\rightarrow 0$ in $L^{p}(\R^N)$.
 Fix $\vartheta=[(1+\eta^2)(1+\bar{c}_{\varepsilon})/(1-\eta)^2\bar{c}_{\varepsilon}]^{1/2}$.
 By virtue of (F1), for $\epsilon =1/4Q_{\max}(\vartheta C_3)^2>0$, there exists a $C_{\epsilon}>0$ such that
 \begin{eqnarray*}
   F(z)\le \epsilon |z|^2+C_{\epsilon}|z|^p,  \ \ \ \ \forall \ z\in \R^2.
 \end{eqnarray*}
 Hence, it follows that
 \begin{eqnarray}\label{F901}
   &     & \limsup_{n\to\infty}\int_{\R^N}Q_{\varepsilon}(x)F(\vartheta {z_n^{\varepsilon}}^{+})\mathrm{d}x\nonumber\\
   & \le & \epsilon Q_{\max}\vartheta^2\limsup_{n\to\infty}\int_{\R^N}|{z_n^{\varepsilon}}^{+}|^2\mathrm{d}x
             +C_{\epsilon}Q_{\max}\vartheta^p\limsup_{n\to\infty}\int_{\R^N}|{z_n^{\varepsilon}}^{+}|^p\mathrm{d}x\nonumber\\
   & \le & \epsilon Q_{\max}(\vartheta C_3)^2=\frac{1}{4}.
 \end{eqnarray}
 From \x{L22}, \x{Ce3}, \x{ab4} and \x{F901}, one has
 \begin{eqnarray*}
   \bar{c}_{\varepsilon}+o(1)
   &  =  & \Phi_{\varepsilon}(z_n^{\varepsilon})\\
   & \ge & \frac{(1-\eta)ab\vartheta^2}{a^2+b^2}\|z_n^{\varepsilon}\|^2
              -\int_{\R^N}Q_{\varepsilon}(x)F(\vartheta{z_n^{\varepsilon}}^{+})\mathrm{d}x\\
   &     & \ \  +\frac{1-\vartheta^2}{2}\langle\Phi_{\varepsilon}'(z_n^{\varepsilon}), z_n^{\varepsilon} \rangle
              +\vartheta^2\langle\Phi_{\varepsilon}'(z_n^{\varepsilon}), {z_n^{\varepsilon}}^{-}\rangle\\
   & \ge & \frac{(1-\eta)ab\vartheta^2}{a^2+b^2}\|{z_n^{\varepsilon}}^{+}\|^2 -\frac{1}{4} +o(1)\\
   & \ge & \frac{(1-\eta)^2\bar{c}_{\varepsilon}\vartheta^2}{1+\eta^2} -\frac{1}{4} +o(1) >  \frac{3}{4}+\bar{c}_{\varepsilon}+o(1).
 \end{eqnarray*}
 This contradiction shows that $\delta_{\varepsilon}>0$.  Going if necessary to a subsequence, we may assume the existence
 of $k_n\in \Z^N$ such that $\int_{B_{1+\sqrt{N}}(k_n)}|{z_n^{\varepsilon}}^{+}|^2dx > \frac{\delta_{\varepsilon}}{2}$. Let $\tilde{z}_n^{\varepsilon}(x)=z_n^{\varepsilon}(x+k_n)$ and $\tilde{\zeta}_n^{\varepsilon}(x)=\zeta_n^{\varepsilon}(x+k_n)$. Then $\|\tilde{z}_n^{\varepsilon}\|=\|z_n^{\varepsilon}\|$ and
 \begin{equation}\label{F903}
   \int_{B_{1+\sqrt{N}}(0)}|(\tilde{z}_n^{\varepsilon})^{+}|^2dx> \frac{\delta_{\varepsilon}}{2}.
 \end{equation}
 Passing to a subsequence, we have $\tilde{z}_n^{\varepsilon}\rightharpoonup \tilde{z}^{\varepsilon}$
 and $(\tilde{z}_n^{\varepsilon})^{+}\rightharpoonup (\tilde{z}^{\varepsilon})^{+}$ in $E$,
 $(\tilde{z}_n^{\varepsilon})^{+}\rightarrow (\tilde{z}^{\varepsilon})^{+}$ in $L^{s}_{\mathrm{loc}}(\R^N)$, $2\le s<2^*$
 and $(\tilde{z}_n^{\varepsilon})^{+}\rightarrow (\tilde{z}^{\varepsilon})^{+}$ a.e. on $\R^N$.
 Obviously, \x{F903} implies that $(\tilde{z}^{\varepsilon})^{+}\ne 0$. Let $\eta_n^{\varepsilon}:=
 \tilde{\zeta}_n^{\varepsilon}/t_n^{\varepsilon}$. Hence, it follows from \x{Phh}, \x{L85}, (F1) and (F3) that
 \begin{eqnarray}\label{F905}
    0 &  =  & \langle \hat{\Phi}'(t_n^{\varepsilon}z_n^{\varepsilon}+\zeta_n^{\varepsilon}),
                t_n^{\varepsilon}z_n^{\varepsilon}+\zeta_n^{\varepsilon} \rangle
                 =\langle \hat{\Phi}'(t_n^{\varepsilon}\tilde{z}_n^{\varepsilon}+\tilde{\zeta}_n^{\varepsilon}),
                t_n^{\varepsilon}\tilde{z}_n^{\varepsilon}+\tilde{\zeta}_n^{\varepsilon} \rangle\nonumber\\
      &  =  & \langle \hat{\Phi}'(t_n^{\varepsilon}(\tilde{z}_n^{\varepsilon}+\eta_n^{\varepsilon})),
                t_n^{\varepsilon}(\tilde{z}_n^{\varepsilon}+\eta_n^{\varepsilon}) \rangle\nonumber\\
      &  =  &  \left[\frac{2ab}{a^2+b^2}\left(\|(\tilde{z}_n^{\varepsilon})^{+}\|^2
                -\|(\tilde{z}_n^{\varepsilon})^{-}+\eta_n^{\varepsilon}\|^2\right)
                +\hat{V}\int_{\R^{N}}|\tilde{z}_n^{\varepsilon}+\eta_n^{\varepsilon}|^2\mathrm{d}x\right](t_n^{\varepsilon})^2\nonumber\\
      &     &   -Q(0)\int_{\R^{N}}F_z(t_n^{\varepsilon}(\tilde{z}_n^{\varepsilon}+\eta_n^{\varepsilon}))\cdot
                  t_n^{\varepsilon}(\tilde{z}_n^{\varepsilon}+\eta_n^{\varepsilon})\mathrm{d}x\nonumber\\
     & \le &  \left[\frac{2ab}{a^2+b^2}\left(\|(\tilde{z}_n^{\varepsilon})^{+}\|^2
                -\|(\tilde{z}_n^{\varepsilon})^{-}+\eta_n^{\varepsilon}\|^2\right)
                +\hat{V}\int_{\R^{N}}|\tilde{z}_n^{\varepsilon}+\eta_n^{\varepsilon}|^2\mathrm{d}x\right](t_n^{\varepsilon})^2\nonumber\\
      &     &   -2Q(0)\int_{\R^{N}}F(t_n^{\varepsilon}(\tilde{z}_n^{\varepsilon}+\eta_n^{\varepsilon}))\mathrm{d}x\nonumber\\
      & \le &  \left[\frac{2ab}{a^2+b^2}\left(\|(\tilde{z}_n^{\varepsilon})^{+}\|^2
                -\|(\tilde{z}_n^{\varepsilon})^{-}+\eta_n^{\varepsilon}\|^2\right)
                +\hat{V}\int_{\R^{N}}|\tilde{z}_n^{\varepsilon}+\eta_n^{\varepsilon}|^2\mathrm{d}x\right](t_n^{\varepsilon})^2\nonumber\\
      &     &   -2\mathcal{C}_0Q(0)(t_n^{\varepsilon})^{\mu}\int_{\R^{N}}F_0(\tilde{z}_n^{\varepsilon}+\eta_n^{\varepsilon})\mathrm{d}x,
 \end{eqnarray}
 which, together with \x{ab3}, implies that
 \begin{eqnarray}\label{F906}
    0 & \le &  \frac{2ab}{a^2+b^2}\left(\|(\tilde{z}_n^{\varepsilon})^{+}\|^2
                -\|(\tilde{z}_n^{\varepsilon})^{-}+\eta_n^{\varepsilon}\|^2\right)
                +\hat{V}\int_{\R^{N}}|\tilde{z}_n^{\varepsilon}+\eta_n^{\varepsilon}|^2\mathrm{d}x\nonumber\\
      & \le &  \frac{2(1+\eta^2)ab}{(1-\eta)(a^2+b^2)}\|(\tilde{z}_n^{\varepsilon})^{+}\|^2
                -\frac{(1-\eta)ab}{a^2+b^2}\|(\tilde{z}_n^{\varepsilon})^{-}+\eta_n^{\varepsilon}\|^2.
 \end{eqnarray}
 This shows that $\{\|\eta_n^{\varepsilon}\|\}_{n=1}^{\infty}$ is bounded in $E^{-}$.  Passing to a subsequence, we have $\eta_n^{\varepsilon}\rightharpoonup \eta^{\varepsilon}$ in $E^{-}$, $\eta_n^{\varepsilon}\rightarrow \eta^{\varepsilon}$ in
 $L^{s}_{\mathrm{loc}}(\R^N)$, $2\le s<2^*$ and $\eta_n^{\varepsilon} \rightarrow \eta^{\varepsilon}$ a.e. on $\R^N$.

 \par
    Since $(\tilde{z}^{\varepsilon}+\eta^{\varepsilon})^{+}=(\tilde{z}^{\varepsilon})^{+}\ne 0$, it follows from (F3)
 that $\int_{\R^N}F_0(\tilde{z}^{\varepsilon}+\eta^{\varepsilon})\mathrm{d}x>0$, which, together with \x{F905}, implies
 that $\{t_n^{\varepsilon}\}_{n=1}^{\infty}$ is bounded, and so $\{\|\zeta_n^{\varepsilon}\|\}_{n=1}^{\infty}$ is also bounded.
 Therefore, there exist $K_1(\varepsilon)>0$ and $K_2(\varepsilon)>0$ such that \x{tnva} holds.
 \qed

 \vskip4mm
 \par\noindent
 {\bf Theorem 4.10.}\ \ {\it Assume that $V$, $Q$ and $F$ satisfy} (V0), (V1) {\it and} (F1)-(F3) . {\it Then
 for $\varepsilon\in (0, \varepsilon_0]$, problem \x{HS6} has a solution $z_{\varepsilon}\in E$ such that
 $\Phi_{\varepsilon}(z_{\varepsilon})=\inf_{\mathcal{N}_{\varepsilon}^{-}}\Phi_{\varepsilon}>0$.}

 \vskip2mm
 \par\noindent
 {\bf Proof.} \  By Lemmas 4.7 and 4.8, there exists a bounded sequence
 $\{z_n^{\varepsilon}\}=\{(u_n^{\varepsilon}, v_n^{\varepsilon})\}\subset E$ satisfying \x{Ce3}. Thus there exists a $z^{\varepsilon}
 =(u^{\varepsilon}, v^{\varepsilon})\in E$ such that $z_n^{\varepsilon}\rightharpoonup z^{\varepsilon}$. Next, we prove that
 $z^{\varepsilon}\ne 0$ for all $\varepsilon\in (0, \varepsilon_0]$.

 \par
   Arguing by contradiction, suppose that $z^{\varepsilon}=0$ for some $\varepsilon\in (0, \varepsilon_0]$,
 i.e. $z^{\varepsilon}_n\rightharpoonup 0$ in $E$, and so $z^{\varepsilon}_n\rightarrow 0$ in $L^{s}_{\mathrm{loc}}(\R^N)$,
 $2\le s<2^*$ and $z^{\varepsilon}_n\rightarrow 0$ a.e. on $\R^N$.

 \par
    We first prove that $\zeta_n^{\varepsilon}\rightharpoonup 0$ in $E^{-}$. Since $\{\|\zeta_n^{\varepsilon}\|\}_{n=1}^{\infty}$
 is bounded, passing to a subsequence we may assume that $\zeta_n^{\varepsilon}\rightharpoonup \zeta^{\varepsilon}$
 in $E^{-}$. By Lemma 4.9, $\{t_n^{\varepsilon}\}$ is bounded, then $t_n^{\varepsilon}z_n^{\varepsilon}+\zeta_n^{\varepsilon}
 \rightharpoonup \zeta^{\varepsilon}$. By Brezis-Lieb's Lemma (see \cite[Lemma 1.32]{Wi}), one can demonstrate that
 \begin{equation}\label{T403}
   \hat{\Phi}(t_n^{\varepsilon}z_n^{\varepsilon}+\zeta_n^{\varepsilon})
     -\hat{\Phi}(t_n^{\varepsilon}z_n^{\varepsilon}+\zeta_n^{\varepsilon}-\zeta^{\varepsilon})-\hat{\Phi}(\zeta^{\varepsilon})=o(1),
 \end{equation}
 which, together with \x{L02}, (V0) and (V1), yields
 \begin{equation}\label{T32}
   o(1) \ge -\hat{\Phi}(\zeta^{\varepsilon})\ge \frac{ab}{a^2+b^2}\|\zeta^{\varepsilon}\|^2
               -\frac{\hat{V}}{2}\int_{\R^{N}}|\zeta^{\varepsilon}|^2\mathrm{d}x
               \ge \frac{(1-\eta)ab}{a^2+b^2}\|\zeta^{\varepsilon}\|^2.
 \end{equation}
 This shows that $\zeta^{\varepsilon}=0$, i.e. $\zeta_n^{\varepsilon}\rightharpoonup 0$ in $E^{-}$.
 By \x{Phh}, \x{L02}, \x{VQR}, \x{Ce3}, \x{L85} and (F1), we have
 \begin{eqnarray}\label{T32}
   {c}_{\varepsilon}+o(1)
    & \ge  & \bar{c}_{\varepsilon}+o(1)=\Phi_{\varepsilon}(z_n^{\varepsilon})\nonumber\\
    & \ge & \Phi_{\varepsilon}(t_n^{\varepsilon}z_n^{\varepsilon}+\zeta_n^{\varepsilon})
             +\frac{1-(t_n^{\varepsilon})^2}{2}\langle\Phi_{\varepsilon}'(z_n^{\varepsilon}), z_n^{\varepsilon}\rangle
             -t_n^{\varepsilon}\langle\Phi_{\varepsilon}'(z_n^{\varepsilon}), \zeta_n^{\varepsilon}\rangle\nonumber\\
    &  =  & \Phi_{\varepsilon}(t_n^{\varepsilon}z_n^{\varepsilon}+\zeta_n^{\varepsilon})+o(1)\nonumber\\
    &  =  & \hat{\Phi}(t_n^{\varepsilon}z_n^{\varepsilon}+\zeta_n^{\varepsilon})
              +\frac{1}{2}\int_{\R^{N}}\left(V_{\varepsilon}(x)-\hat{V}\right)|t_n^{\varepsilon}z_n^{\varepsilon}
              +\zeta_n^{\varepsilon}|^2\mathrm{d}x\nonumber\\
    &     &   +\int_{\R^{N}}\left[Q(0)-Q_{\varepsilon}(x)\right]F(t_n^{\varepsilon}z_n^{\varepsilon}
              +\zeta_n^{\varepsilon})\mathrm{d}x+o(1)\nonumber\\
    & \ge & \hat{c}+\frac{1}{2}\int_{|x|\le R_0/\varepsilon}\left(V_{\varepsilon}(x)-\hat{V}\right)|t_n^{\varepsilon}z_n^{\varepsilon}
              +\zeta_n^{\varepsilon}|^2\mathrm{d}x\nonumber\\
    &     &   +\int_{|x|\le R_0/\varepsilon}\left[Q(0)-Q_{\varepsilon}(x)\right]F(t_n^{\varepsilon}z_n^{\varepsilon}
              +\zeta_n^{\varepsilon})\mathrm{d}x+o(1)\nonumber\\
    & \ge & \hat{c}-\frac{2V_{\max}-V_{\infty}-V(0)}{4}\int_{|x|\le R_0/\varepsilon}|t_n^{\varepsilon}z_n^{\varepsilon}
              +\zeta_n^{\varepsilon}|^2\mathrm{d}x\nonumber\\
    &     &   +\left[Q(0)-Q_{\max}\right]\int_{|x|\le R_0/\varepsilon}F(t_n^{\varepsilon}z_n^{\varepsilon}
              +\zeta_n^{\varepsilon})\mathrm{d}x+o(1)\nonumber\\
    &  =  & \hat{c}+o(1).
 \end{eqnarray}
 On the other hand, by Lemmas 4.3 and 4.6, one has
 $$
   \hat{c}\ge c_0+\hat{\delta}\ge c_{\varepsilon}+\frac{1}{4}\hat{\delta},
 $$
 which contradicts to \x{T32}. Therefore, $z^{\varepsilon}\ne 0$ for all $\varepsilon\in (0, \varepsilon_0]$.
 In a standard way, we can certify that $\Phi_{\varepsilon}'(z^{\varepsilon})=0$ and
 $\Phi_{\varepsilon}(z^{\varepsilon})=c_{\varepsilon}=\inf_{\mathcal{N}_{\varepsilon}^{-}}\Phi_{\varepsilon}$. This shows that
 $z^{\varepsilon}\in E$ is a solution for problem \x{HS6} with $\Phi_{\varepsilon}(z^{\varepsilon})=\inf_{\mathcal{N}_{\varepsilon}^{-}}\Phi_{\varepsilon}>0$.
 \qed

 \vskip4mm
 \par\noindent
 {\bf Proof of Theorem 1.6.} \  For $\varepsilon\in (0, \varepsilon_0]$, Theorem 4.10 implies that problem \x{HS6}
 has a solution $z_{\varepsilon}\in E$ such that $\Phi_{\varepsilon}(z_{\varepsilon})=\inf_{\mathcal{N}_{\varepsilon}^{-}}
 \Phi_{\varepsilon}>0$. Then
 $$
   \hat{z}_{\varepsilon}(x)=(\hat{u}_{\varepsilon}(x), \hat{v}_{\varepsilon}(x))
   :=z_{\varepsilon}(\varepsilon^{-1}(x-x_v))=(u_{\varepsilon}(\varepsilon^{-1}(x-x_v)), v_{\varepsilon}(\varepsilon^{-1}(x-x_v)))
 $$
 is a nontrivial solution of problem \x{PHS}.
 \qed

 {}

\end{document}